\documentclass[a4paper, 11pt]{article}

\usepackage{mathrsfs}
\usepackage{epsfig}
\usepackage{amsmath}
\usepackage{amssymb}
\usepackage{latexsym}
\usepackage{amsfonts}
\usepackage{amsthm}
\usepackage[numbers,sort&compress]{natbib}
\usepackage{graphicx}
\ifx\pdfoutput\undefined  \else
 \fi
\usepackage[centerlast]{caption2}
\usepackage{color}

\newcommand {\cross } {\nu_D }

\newtheorem{definition}{Definition}

\newtheorem{lemma}{Lemma}
\newtheorem{theorem}{Theorem}

\numberwithin{figure}{section} \numberwithin{definition}{section}
\numberwithin{observation}{section} \numberwithin{lemma}{section}
\numberwithin{theorem}{section}

\begin{document}

\title{
{Disproof of a conjecture by Erd\H{o}s and Guy on the crossing number
of hypercubes}
\author{
Yuansheng Yang$^1$\thanks{E-mail:
yangys@dlut.edu.cn}, \ \ Guoqing Wang$^2$\thanks{Corresponding E-mail:
gqwang1979@aliyun.com}, \ \ Haoli
Wang$^3$\thanks{E-mail: bjpeuwanghaoli@163.com}, \ \  Yan Zhou$^1$\\
\\
$^1$Department of Computer Science \\
Dalian University of Technology, Dalian, 116024, P. R. China\\
\\
$^2$School of Mathematical Sciences \\
Tiangong University, Tianjin, 300387, P. R. China\\
\\
$^{3}$ College of Computer and Information Engineering \\
Tianjin Normal University, Tianjin, 300387, P. R. China\\
}}

\date{}
\maketitle
\begin{abstract} Let $Q_n$ be the $n$-dimensional hypercube, and let
${\rm cr}(Q_n)$ be the \textit{crossing number} of $Q_n$. Erd\H{o}s and Guy in 1973
conjectured
the following equality: ${\rm
cr}(Q_n)=\frac{5}{32}4^n-\lfloor\frac{n^2+1}{2}\rfloor 2^{n-2}$. In
this paper, we construct a drawing of $Q_n$ with less crossings when $n>6$,
which implies that for $n>6$ we have a strict inequality.
\end{abstract}

\medskip


\section{Introduction}

Let $G$ be a simple connected graph with vertex set $V(G)$ and edge
set $E(G)$. The  \textit{crossing number} of a graph $G$, denoted
${\rm cr}(G)$, is the minimum possible number of edge crossings in a
drawing of $G$ in the plane. The notion of crossing number is a
central one for Topological Graph Theory and has been studied
extensively  by mathematicians including Erd\H{o}s, Guy, Harary,
Tur\'{a}n and Tutte, et al. (see
\cite{EG73,Harary,SSSV,Turan77,Tutte70}).

The investigation on the crossing number of graphs is a difficult problem. In 1973, Erd\H{o}s and Guy \cite{EG73}
wrote, ``{\sl Almost all questions that one can ask about crossing
numbers remain unsolved.}''  Garey and Johnson \cite{GJ83}
proved that computing the crossing number is NP-complete. In this
field, Eggleton and  Guy \cite{Egg70} in 1970 announced a drawing
of the $n$-dim hypercube $Q_n$ with exactly $\frac{5}{32}4^n-\lfloor\frac{n^2+1}{2}\rfloor
2^{n-2}$ crossings, which implies \begin{equation}\label{equation
conjeture of Erdos and Guy} {\rm
cr}(Q_n)\leqq\frac{5}{32}4^n-\lfloor\frac{n^2+1}{2}\rfloor 2^{n-2},
\end{equation}
where $Q_n$ denotes the $n$-dimensional hypercube.

Not long afterward, in 1971 a gap was found  \cite{Guytrend} in the drawing given by Eggleton and Guy.
We still quote: \emph{``but a gap has been found in the description
of the construction, so this must also remain a conjecture. We again
conjecture equality in \eqref{equation conjeture of Erdos and
Guy}.''} (P. Erd\H{o}s and R.K. Guy \cite{EG73})

For a long time, the \emph{equality} conjectured by
Erd\H{o}s and Guy was believed to be true, even though a proof appeared out of reach.
The past results on ${\rm
cr}(Q_n)$, including ${\rm cr}(Q_3)=0$ (trivial), ${\rm cr}(Q_4)=8$
(see \cite{DeaRic95}), and the obtained best drawings for $Q_n$ with
$n=5,6,7,8$ (see \cite{FF00}, \cite{M91}) seem to support \emph{the
equality}.

Very recently in 2008, Faria, Figueiredo,
S\'{y}kora and Vr\v{t}o \cite{FFSV08} announced a drawing for which
the number of crossings coincides with
$\frac{5}{32}4^n-\lfloor\frac{n^2+1}{2}\rfloor 2^{n-2}$, giving further support to the conjecture.

In this paper, we
construct drawings of the hypercubes $Q_n$ with, for $n>6$, fewer crossings, disproving the conjecture.
We prove the following upper bound for the crossing
number of hypercubes.
\begin{theorem}\label{Theorem main theorem}
${\rm cr}(Q_n)\leq \frac{139}{896} 4^{n}-\lfloor\frac{n^2+1}{2}\rfloor 2^{n-2}+\frac{4}{7}\cdot 2^{3\lfloor\frac{n}{2}\rfloor -n}+\lambda_n(\frac{4^n}{24576} -\frac{4^{\lfloor\frac{ n}{2}\rfloor}}{6}),$ where $$\begin{array}{llll} \lambda_n= &\left
\{\begin{array}{llll}
               0, & \mbox{ if } 5\leq n\leq 12; \\
               1, & \mbox{ if } n\geq 13. \\
                       \end{array}
           \right. \\
\end{array}$$
\end{theorem}

Our proposed family of drawings for $Q_n,$ when $n=5$ and $n=6$, has the exact number of crossings conjectured by Erd\H{o}s and Guy, and when $n>6$ has less crossing number than that conjectured.  In the concluding section of this paper, we will compare the asymptotic behavior of our proposed new function for the upper bound of ${\rm cr}(Q_n)$ to previous function.

\section{Notation}

A drawing of $G$ is said to be a {\it good} drawing, provided that
no edge crosses itself, no adjacent edges cross each other, no two
edges cross more than once, and no three edges cross in a point. It
is well known that the crossing number of a graph is attained only
in {\it good} drawings of the graph. So, we always assume that all
drawings throughout this paper are good drawings. Let $D$ be a good
drawing of the graph $G$, and let $A$ and $B$ be two disjoint
subsets of $E(G)$. In the drawing $D$, the number of the crossings
formed by the edges of $A$ with the edges of $B$ is denoted by
$\cross(A,B)$. The number of the crossings between the edges of
$A$ is denoted by $\cross (A)$. In what follows,
$\nu_D(E(G))$ is abbreviated to $\nu(D)$ when it is unambiguous. Let
$u$ be a vertex of $G$, and let $U$ be a vertex subset of $V(G)$. We
define $\mathcal {I}(u)$ to be the edge subset of $E(G)$ consisting of all edges incident with $u$. Let
$$\mathcal{I}(U)=\bigcup\limits_{u\in U}\mathcal {I}(u),$$ and let
$$\partial(U)=\mathcal{I}(U)\setminus E(U),$$
where $E(U)=\{uv\in E(G):u,v\in U\}$.

The $n$-dimensional hypercube $Q_n$ is a graph with the vertex set
$V(Q_n)=\{d_1d_2 \cdots d_n: d_i\in \{0,1\}, i=1,2, \ldots, n\}$,
for which any two vertices $a=a_1a_2\cdots a_n$ and $b=b_1b_2\cdots
b_n$ are adjacent if and only if there exists a unique $i\in \{1,2,\ldots,n\}$
such that $a_i\neq b_i$.
In particular, if the unique $i\in \{1,2,\ldots,n\}$
with $a_i\neq b_i$ is equal to $n-1$, i.e., $a=a_1a_2\cdots a_{n-2}a_{n-1}a_n$ and $b=a_1a_2\cdots a_{n-2}\overline{a_{n-1}}a_n$, we denote $$b=\widehat{a},$$
and
conversely, $$a=\widehat{b}.$$
For any vertex $a=a_1a_2\cdots a_n\in V(Q_n)$ and any binary string $x_1x_2\cdots
x_t$ of length $t$, we define $$a^{(x_1x_2\cdots x_t)}=a_1a_2\cdots a_nx_1x_2\cdots
x_t$$ to be the vertex of $V(Q_{n+t})$.

\section{Proof of Theorem \ref{Theorem main theorem}}

In this section, we prove Theorem \ref{Theorem main theorem} by
constructing a drawing, denoted $\Gamma_n$, of $Q_n$ with the
desired crossings for every integer $n\geq 5$. The constructions of
the drawing $\Gamma_n$ are different according to the parity of $n$.  Hence, we shall introduce the constructions of $\Gamma_n$
in Subsection 3.1 and Subsection 3.2. In Subsection 3.3,
we verify the constructed drawing $\Gamma_n$ has the
desired number of crossings.

\subsection{Construction of the drawing $\Gamma_n$ for all odd $n$}

\noindent $\bullet$ {\sl Throughout this subsection, we use $n$ as
an odd integer no less than $5$.}

\begin{figure}[!ht]
\centering
\includegraphics[scale=1.0]{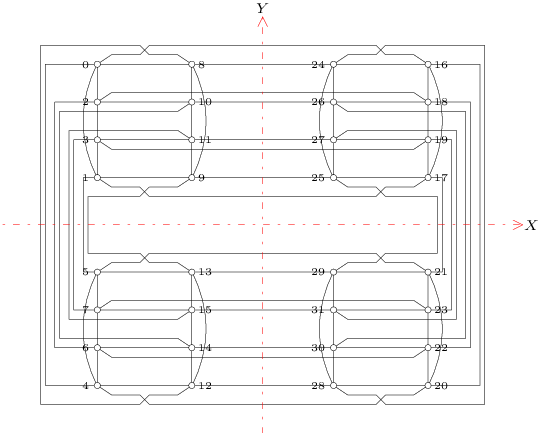}
\caption{\small{The drawing $\Gamma_5$ with 56
crossings}}\label{fig:gamma 5}
\end{figure}

The desired drawing $\Gamma_n$ of $Q_n$ will be constructed recursively (see Figure \ref{fig:gamma 5} for the drawing $\Gamma_5$). Note that
the drawing $\Gamma_5$ for $Q_5$ we give here is different to the previous drawing for $Q_5$ in \cite{FF00,M91}.  The drawing in Figure \ref{fig:gamma 5} has the same number of crossings as the drawing obtained in \cite{FF00,M91}, however, it can help us to inductively construct our general drawings for $Q_n$ with fewer crossings than the crossing number conjectured by Erd\H{o}s and Guy.

The idea of constructing the drawing $\Gamma_{n+2}$  from the obtained drawing $\Gamma_n$ is as follows. By preprocessing the drawing $\Gamma_n$, we get a new drawing of $Q_n$, denoted $\Gamma_n^*$, with some properties which are helpful for subsequent procedures.  Replacing each vertex $u$ in $\Gamma_n^*$ by four new vertices $u^{(00)},u^{(10)},u^{(11)},u^{(01)}\in V(Q_{n+2})$ in the `very small neighborhood' of the original location of $u$, and replacing each edge $uv$ in $\Gamma_n^*$ by a `bunch' (in this paper, a bunch means edges drawn in parallel close to each other) of four new edges $u^{(00)}v^{(00)}, u^{(10)}v^{(10)}, u^{(11)}v^{(11)}, u^{(01)}v^{(01)}\in E(Q_{n+2})$ drawn along the original route of $uv$, we get a {\bf transitional} drawing of $Q_{n+2}$, denoted $\widetilde{\Gamma}_{n+2}$. Next, by modifying the `routes' of some edges in $\widetilde{\Gamma}_{n+2}$, we decrease the crossings and obtain the desired drawing $\Gamma_{n+2}$ of $Q_{n+2}$. In what follows,  we shall give a detailed description of the construction.

In order to accurately describe the process, all the drawings in this paper will be given in the $2$-dimensional Euclidean plane
$\mathbb{R}\times \mathbb{R}$. For any vertex $u$ in some drawing, by $X_u$ and $Y_u$ we denote the $X$
and $Y$-coordinates of $u$ in $\mathbb{R}\times \mathbb{R}$. Three different drawings $\Gamma_n,\Gamma_n^*$ and $\widetilde{\Gamma}_{n}$ ($\widetilde{\Gamma}_{n}$ exists only for $n\geq 7$) of $Q_n$ will appear later, for convenience, we use the notion $\mathcal{D}_n$ to denote any one of the above three possible drawings $\Gamma_n,\Gamma_n^*$ and $\widetilde{\Gamma}_{n}$.

In this paper, any drawing $\mathcal{D}_n$ shares the following inductive rule for the arrangements of vertices:

For $n=5$, the locations of vertices in $\mathcal{D}_5$ (including $\Gamma_5$ and $\Gamma_5^*$) are depicted
in Figure \ref{fig:gamma 5}, in particular, $$\{(X_u,Y_u):u
\mbox{ is a vertex in } \mathcal{D}_5\}=\{-2,-1,1,2\}\times
\{-4,-3,-2,-1,1,2,3,4\}.$$

\noindent $\bullet$ {\small For convenience, in the rest of
this paper, any vertex $a=a_1a_2\cdots a_n\in V(Q_n)$ drawn in figures
will be represented by the corresponding decimal number
$2^{n-1}a_1+2^{n-2}a_2+\cdots +2^0a_n$.}

\noindent $\bullet$ {\small Throughout this paper, we always let $\mathcal{N}$ denote} some fixed large positive integer.

\begin{figure}[!ht]
\centering
\includegraphics[scale=0.85]{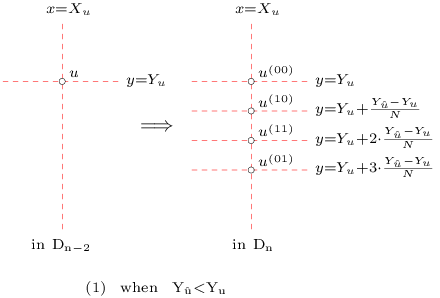}\hspace{50pt}
\includegraphics[scale=0.85]{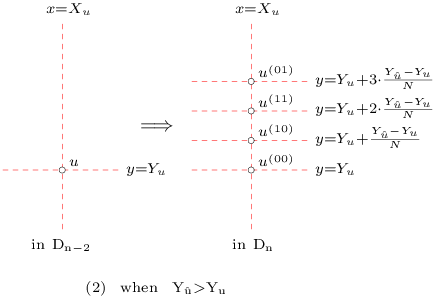}
\caption{\small{Inductive rule for the arrangement of vertices}}\label{fig:inductiveruleforvertices}
\end{figure}

Suppose $n>5$. Take an
arbitrary vertex $u$ in $\mathcal{D}_{n-2}$.  The four vertices extended from $u$, say
$u^{(00)},u^{(10)},u^{(11)},u^{(01)}$, in $\mathcal{D}_n$ will be
located in $\mathbb{R}\times \mathbb{R}$ such that
\begin{equation}\label{equation X for odd n}
X_{u^{(00)}}=X_{u^{(10)}}=X_{u^{(11)}}=X_{u^{(01)}}=X_{u}
\end{equation} and
\begin{equation}\label{equation Y for odd n}
\begin{array}{llll} &\left \{\begin{array}{llll}
               Y_{u^{(00)}}=Y_{u}; \\
               Y_{u^{(10)}}=Y_{u}+\frac{Y_{\widehat{u}}-Y_{u}}{\mathcal {N}};\\
               Y_{u^{(11)}}=Y_{u}+2\cdot
               \frac{Y_{\widehat{u}}-Y_{u}}{\mathcal{N}};\\
               Y_{u^{(01)}}=Y_{u}+3\cdot
               \frac{Y_{\widehat{u}}-Y_{u}}{\mathcal{N}}.\\
                       \end{array}
           \right. \\
\end{array}
\end{equation}
 In other words, $u^{(00)},u^{(10)},u^{(11)},u^{(01)}$ in $\mathcal{D}_n$ are drawn at the line where $u$ lies, and are located at the `very small neighborhood' of the original location of $u$ in $\mathcal{D}_{n-2}$, noticing that $\mathcal{N}$ is some fixed large positive integer. To make it clear, we illustrate the above inductive rule in Figure \ref{fig:inductiveruleforvertices}.

{\small In some cases, we shall denote the four vertices $u^{(00)},u^{(10)},u^{(11)},u^{(01)}$ in Figure  \ref{fig:inductiveruleforvertices} to be $\mathcal{P}_1(u)$, $\mathcal{P}_2(u)$, $\mathcal{P}_3(u)$, $\mathcal{P}_4(u)$ where $Y_{\mathcal{P}_1(u)}>  Y_{\mathcal{P}_2(u)}>  Y_{\mathcal{P}_3(u)}> Y_{\mathcal{P}_4(u)}$, i.e., $$(\mathcal{P}_1(u),  \mathcal{P}_2(u),  \mathcal{P}_3(u), \mathcal{P}_4(u))=(u^{(00)},u^{(10)},u^{(11)},u^{(01)})$$ for Diagram (1), and $$(\mathcal{P}_1(u),  \mathcal{P}_2(u),  \mathcal{P}_3(u), \mathcal{P}_4(u))=(u^{(01)},u^{(11)},u^{(10)},u^{(00)})$$ for Diagram (2), respectively.}

By the above inductive rule for the
arrangements of vertices, we have the following.

\noindent \textbf{Claim A.} \ {\sl Let $u_1$ and $u_2$ be two adjacent vertices in $\mathcal{D}_n$.
Then,

{\rm (i)}. \  either $X_{u_1}=X_{u_2}$ or $Y_{u_1}=Y_{u_2}$ holds;

{\rm (ii)}. \ if $u_1=\widehat{u_2}$, then $X_{u_1}=X_{u_2}$ and both vertices $u_1,u_2$ are drawn next to each other at the line $x=X_{u_1}$;

{\rm (iii)}. \  $(Y_{\widehat{u_1}}-Y_{u_1})\cdot (Y_{\widehat{u_2}}-Y_{u_2})$ is negative or  positive according to $X_{u_1}=X_{u_2}$ or $Y_{u_1}=Y_{u_2}$.}

\medskip

\noindent {\bf Proof of Claim A.} \  We first prove that  Conclusion {\rm (ii)} holds.
Say $u_1=\widehat{u_2}$. For $n=5$, Conclusion {\rm (ii)}
holds by Figure \ref{fig:gamma 5}. Hence, we need only to consider the case that $n>5$, i.e., $n\geq 7$ by recalling that $n$ is odd in this subsection.
Then there exists some vertex $v$ in $\mathcal{D}_{n-2}$ such that $$\mbox{ either }\ \ \{u_1,u_2\}=\{v^{(00)},v^{(10)}\} \ \ \mbox{ or }\ \ \{u_1,u_2\}=\{v^{(11)},v^{(01)}\}.$$
From Figure \ref{fig:inductiveruleforvertices}, we see immediately that $u_1$ and $u_2$ are drawn next to each other at the same line $x=X_{u_1}=X_{u_2}$ in $\mathcal{D}_n$, which proves Conclusion {\rm (ii)}.

Next we shall prove Conclusion {\rm (i)} and Conclusion {\rm (iii)} by induction on $n$. Similarly as above, both conclusions for the case of $n=5$ can be verified in Figure \ref{fig:gamma 5}. We assume $n\geq 7$.
Let
$$u_1=v_1^{(a_1a_2)}$$
and
$$u_2=v_2^{(b_1b_2)}$$
where $v_1,v_2\in V(Q_{n-2})$ and $a_1a_2,b_1b_2\in \{00,10,11,01\}.$
Since $u_1,u_2$ are adjacent, we have that
$$v_1=v_2 \ \ \ \mbox{or} \ \ \  a_1a_2=b_1b_2.$$ Then we distinguish two cases.

\textbf{Case 1.} \ $v_1=v_2$.

By Conclusion {\rm (ii)}, we may assume without loss of generality that $u_1\neq \widehat{u_2}$, and thus,
$$\mbox{ either }\ \ \{u_1,u_2\}=\{v_1^{(10)},v_1^{(11)}\} \ \ \mbox{ or }\ \ \{u_1,u_2\}=\{v_1^{(01)},v_1^{(00)}\}.$$
From Figure \ref{fig:inductiveruleforvertices},
we see that
$X_{u_1}=X_{u_2}=X_{v_1}$ and $(Y_{\widehat{u_1}}-Y_{u_1})\cdot (Y_{\widehat{u_2}}-Y_{u_2})<0$, and that Conclusion {\rm (i)} and Conclusion {\rm (iii)} follow.

\textbf{Case 2.} \ $a_1a_2=b_1b_2$.

 Since $u_1$ and $u_2$ are adjacent, it follows that $v_1$ and $v_2$ are adjacent vertices of $V(Q_{n-2})$. By applying the induction hypothesis for Conclusion (i) and (iii) on $n-2$,
we have that in $\mathcal{D}_{n-2}$,
either
\begin{equation}\label{equation <0}
X_{v_1}=X_{v_2}\ \ \ \mbox{ and } \ \ \ (Y_{\widehat{v_1}}-Y_{v_1})\cdot (Y_{\widehat{v_2}}-Y_{v_2})<0,
\end{equation}
 or
 \begin{equation}\label{equation >0}
 Y_{v_1}=Y_{v_2} \ \ \ \mbox{ and } \ \ \ (Y_{\widehat{v_1}}-Y_{v_1})\cdot (Y_{\widehat{v_2}}-Y_{v_2})>0.
 \end{equation}
By Conclusion {\rm (ii)},  $X_{\widehat{v_i}}=X_{v_i}$ and that $\widehat{v_i}$ and $v_i$ are drawn next to each other in the line $x=X_{v_i}$ where $i=1,2$.

\begin{figure}[!ht]
\centering
\includegraphics[scale=0.9]{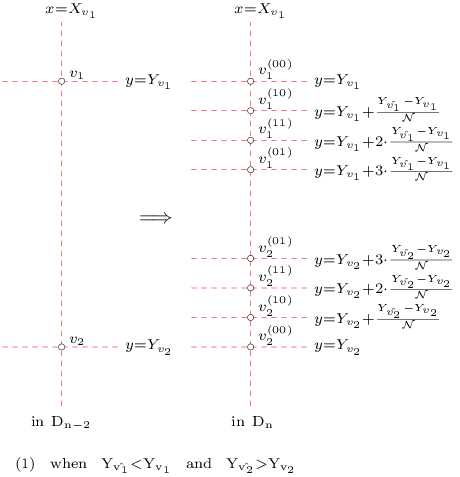}\hspace{40pt}
\includegraphics[scale=0.9]{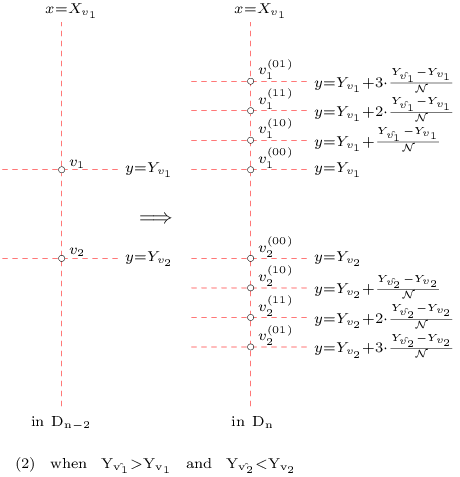}
\caption{\small{Auxiliary diagrams illustrating Claim A}}
\label{fig:locationofadjacentvertices}
\end{figure}

Suppose \eqref{equation <0} holds.  We conclude that $X_{u_1}=X_{u_2}$ and $(Y_{\widehat{u_1}}-Y_{u_1})\cdot (Y_{\widehat{u_2}}-Y_{u_2})<0$ (see Figure \ref{fig:locationofadjacentvertices}).

\begin{figure}[!ht]
\centering
\includegraphics[scale=0.75]{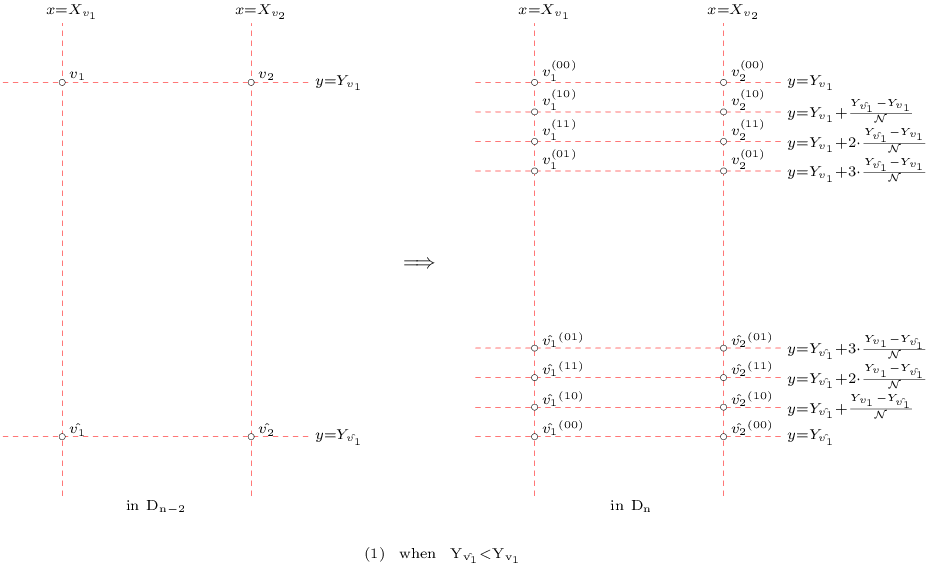}

\vspace{20pt}

\includegraphics[scale=0.75]{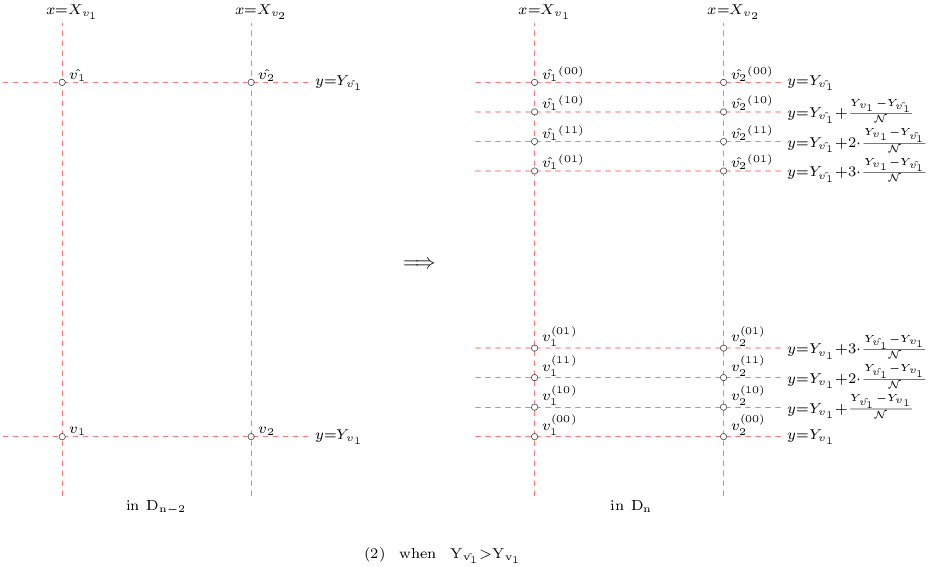}
\caption{\small{Auxiliary diagrams illustrating Claim A}}
\label{fig:Ralativelocationadvertices}
\end{figure}

Suppose \eqref{equation >0} holds. Since $v_1$ and $v_2$ are adjacent vertices of $V(Q_{n-2})$, we have that $\widehat{v_1}$ and $\widehat{v_2}$ are adjacent too. Since $X_{v_1}\neq X_{v_2}$, it follows from Conclusion {\rm (ii)} that $X_{\widehat{v_1}}\neq X_{\widehat{v_2}}$. By applying the induction hypothesis for Conclusion (i) on $n-2$, we have that $Y_{\widehat{v_1}}=Y_{\widehat{v_2}}$. Then we conclude that $Y_{u_1}=Y_{u_2}$ and $(Y_{\widehat{u_1}}-Y_{u_1})\cdot (Y_{\widehat{u_2}}-Y_{u_2})>0$ (see Figure \ref{fig:Ralativelocationadvertices}).
\textbf{This proves Claim A.} \qed

\bigskip

Then, we shall characterize the drawing of the edges.
As explained in the beginning of this subsection, for any odd number
$n\geq 5$, the drawing of edges in $\Gamma_{n+2}$ will be constructed from
$\Gamma_{n}$ inductively.
We shall obtain the drawing of edges in $\Gamma_{n+2}$ from $\Gamma_{n}$ in three
steps.  We give a sketch of the three steps here, and give the accurate description later.

In the first step, we obtain another drawing of $Q_n$,
denoted $\Gamma_n^*$, for which in the `very small neighborhood' of any vertex $u$ of $V(Q_n)$, the number of edges incident with $u$ which are drawn on the left of the line $x=X_u$ is almost the same as the number of edges incident with $u$ which are drawn on the right of the line $x=X_u$.

In the second step, each vertex $u$ in $\Gamma_n^*$ will be replaced by the  $4$-cycle $u^{(00)}u^{(10)}u^{(11)}u^{(01)}$ drawn in the `very small neighborhood' of the original location of $u$ and `precisely' at the line $x=X_u$, and every edge $uv$ in $\Gamma_n^*$
will be replaced by a `bunch' of four new edges
$u^{(00)}v^{(00)}$, $u^{(10)}v^{(10)}$, $u^{(11)}v^{(11)}$, $u^{(01)}v^{(01)}$ drawn along the original route of $uv$. This gives the {\sl transitional} drawing $\widetilde{\Gamma}_{n+2}$ of $Q_{n+2}$.

In the final step, we decrease the crossings in the drawing $\widetilde{\Gamma}_{n+2}$ by adjusting the route of some edges `locally' or `globally' and get the final desired drawing $\Gamma_{n+2}$.  Say $e$ is an edge of  $\{u^{(00)}v^{(00)},u^{(10)}v^{(10)},u^{(11)}v^{(11)},u^{(01)}v^{(01)}\}$ to be adjusted in $\widetilde{\Gamma}_{n+2}$.  By saying to adjust the edge $e$ `locally', we mean altering the drawing of the two parts of the edge $e$ which are located in the `very small neighborhoods' of both two 4-cycles $u^{(00)}u^{(10)}u^{(11)}u^{(01)}$ and $v^{(00)}v^{(10)}v^{(11)}v^{(01)}$, but drawing the edge  $e$ still in its original bunch  $\{u^{(00)}v^{(00)},u^{(10)}v^{(10)},u^{(11)}v^{(11)},u^{(01)}v^{(01)}\}$ outside both neighborhoods.
By saying to adjust the edge $e$ `globally', we mean altering the whole route of the edge $e$ and not necessarily drawing the edge $e$ still in its original bunch.

To proceed, we shall need some technical definitions.

\begin{definition} Let $u$ be a vertex and $e$ an edge incident with $u$ in $\mathcal{D}_n$. We call the edge $e$

$\bullet$ \verb"a left arc" with respect to $u$, if the `starting part' of $e$
within the `very small neighborhood' of $u$ is drawn on the left of
the line $x=X_u$;

$\bullet$  \verb"a right arc" with respect to $u$, if the `starting part' of
$e$ within the `very small neighborhood' of $u$ is drawn on the right
of the line $x=X_u$;

$\bullet$ \verb"a below arc" with respect to $u$, if the `starting part' of
$e$ within the `very small neighborhood' of $u$ is drawn below
the line $y=Y_u$;

$\bullet$ \verb"an above arc" with respect to $u$, if the `starting part' of
$e$ within the `very small neighborhood' of $u$ is drawn above
the line $y=Y_u$.

Let $\mathcal{L}_{\mathcal{D}_n}(u), \mathcal{R}_{\mathcal{D}_n}(u), \mathcal{B}_{\mathcal{D}_n}(u), \mathcal{A}_{\mathcal{D}_n}(u),$ be the set of edges incident with $u$ which are {\bf left, right, below, above} arcs with respect to $u$ in the drawing $\mathcal{D}_n$, respectively.
\end{definition}

\begin{figure}[!ht]
\centering
\includegraphics[scale=1.0]{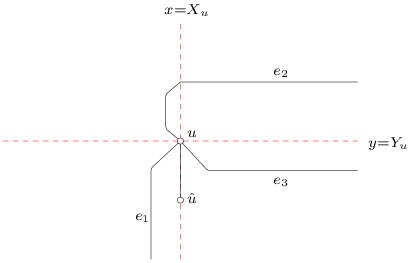}
\caption{\small{The natures of any edge with respect to its ends in $\mathcal{D}_n$}}
\label{fig:Natureofedge}
\end{figure}

Now let's give an example to illustrate the above points.
It can be seen in Figure \ref{fig:Natureofedge} that, $e_1\in \mathcal{L}_{\mathcal{D}_n}(u)\cap  \mathcal{B}_{\mathcal{D}_n}(u)$, i.e., the edge $e_1$ is both a left arc and a below arc with respect to $u$ in the drawing $\mathcal{D}_n$;
$e_2\in \mathcal{L}_{\mathcal{D}_n}(u)\cap  \mathcal{A}_{\mathcal{D}_n}(u)$;  $e_3\in \mathcal{R}_{\mathcal{D}_n}(u)\cap  \mathcal{B}_{\mathcal{D}_n}(u)$; While the edge $u\widehat{u}\in \mathcal{B}_{\mathcal{D}_n}(u)$ and $u\widehat{u}\notin \mathcal{L}_{\mathcal{D}_n}(u)\cup \mathcal{R}_{\mathcal{D}_n}(u)$, i.e., the edge $u\widehat{u}$ is just a below arc with respect to $u$, and is neither a
left arc nor a right arc with respect to $u$ since the edge is drawn precisely along the line $x=X_u$.  It is worth while to note that the above {\bf nature} of any edge $e$ with respect to any of its ends $u$ is determined just by the route of the part of $e$ which is located in the `very small neighborhood' of $u$. That is also the reason why the edge $e_2$ is a left arc rather than a right arc with respect to $u$. It is easy to check that in Figure \ref{fig:Natureofedge}, $|\mathcal{L}_{\mathcal{D}_n}(u)|=2, |\mathcal{R}_{\mathcal{D}_n}(u)|=1, |\mathcal{B}_{\mathcal{D}_n}(u)|=3, |\mathcal{A}_{\mathcal{D}_n}(u)|=1$.

\begin{definition}\label{definition edge self-symmetric} Let
$e=u_1u_2$ be an edge in $\mathcal{D}_n$. We say that the edge
$e$ is \verb"self-symmetric" provided that the following
conditions hold:

If $X_{u_1}=X_{u_2}$, then the edge $e$ is drawn symmetrically with
respect to the line $y=\frac{Y_{u_1}+Y_{u_2}}{2}$, i.e., the part drawn above is a reflection of the part drawn below with the `mirror' $y=\frac{Y_{u_1}+Y_{u_2}}{2}$, in
particular, $e$ is a left (right) arc with respect to $u_1$ if and only if $e$ is a left (right) arc with respect to $u_2$, and $e$ is an above (below) arc with respect to $u_1$ if and only if $e$ is a below (above) arc with respect to $u_2$ (see Figure \ref{fig:Selfsymmx=x}).
\begin{figure}[!ht]
\centering
\includegraphics[scale=0.9]{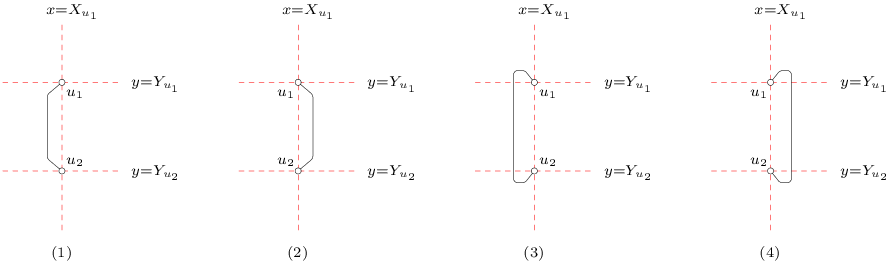}
\caption{\small{The edge $u_1u_2$ is self-symmetric with $X_{u_1}=X_{u_2}$}}\label{fig:Selfsymmx=x}
\end{figure}

If $Y_{u_1}=Y_{u_2}$, then the edge $e$ is drawn symmetrically with
respect to the line $x=\frac{X_{u_1}+X_{u_2}}{2}$, i.e., the part drawn on the left is a reflection of the part drawn on the
right  with the `mirror' $x=\frac{X_{u_1}+X_{u_2}}{2}$, in
particular,
$e$ is a left (right) arc with respect to $u_1$ if and only if $e$ is a right (left) arc with respect to $u_2$, and $e$ is an above (below) arc with respect to $u_1$ if and only if $e$ is an above (below) arc with respect to $u_2$ (see Figure \ref{fig:Selfsymmy=y}).
\begin{figure}[!ht]
\centering
\includegraphics[scale=0.8]{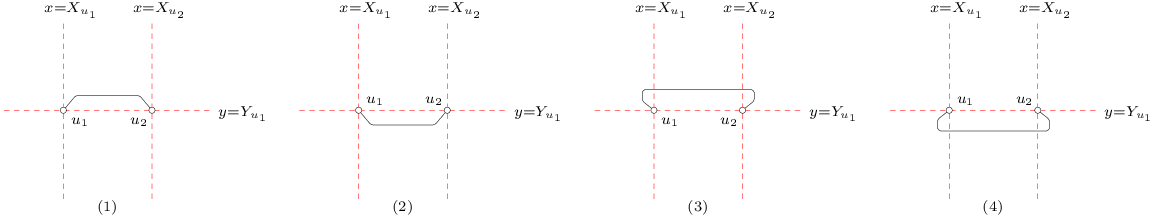}
\caption{\small{The edge $u_1u_2$ is self-symmetric with $Y_{u_1}=Y_{u_2}$}}
\label{fig:Selfsymmy=y}
\end{figure}
\end{definition}

It can be seen that in the drawing $\Gamma_5$ given by Figure \ref{fig:gamma 5}, every edge is self-symmetric. Furthermore, every edge is self-symmetric in any drawing $\mathcal{D}_n$ of $Q_n$ for all $n\geq 5$ given in this paper. {\bf One basic principle} in constructing $\Gamma_{n+2}$ from $\Gamma_n$ is to keep the adjusted edges self-symmetric after any adjustment in each one of the three steps, from $\Gamma_n$ to $\Gamma_n^*$, or from $\Gamma_n^*$ to $\widetilde{\Gamma}_{n+2}$, or from $\widetilde{\Gamma}_{n+2}$ to $\Gamma_{n+2}$. In fact, any adjustment on some edge will be done with respect to its both ends synchronously in order to keep the adjusted edge self-symmetric.
The main reason to keep edges self-symmetric is to make the four new  edges $u^{(00)}v^{(00)},u^{(10)}v^{(10)},u^{(11)}v^{(11)},u^{(01)}v^{(01)}$ do not cross each other for any edge $uv$ in $\Gamma_n$ during the process of constructing $\Gamma_{n+2}$ from $\Gamma_n$, which shall be illustrated more specifically later.

\medskip

Now we give a structure which is crucial in this paper to construct
the drawing $\Gamma_n$ with the desired number of crossings.

\begin{definition} Let $U=\{u_1,u_2,\ldots,u_8\}$ be a set
of eight vertices  in $\mathcal{D}_n$.
 We say that the drawing of
the edge set $\mathcal{I}(U)$ forms \verb"a fundamental structure",
 provided that the drawing of
$\mathcal{I}(U)$ is given as either Diagram (1) or Diagram (2) in
Figure \ref{fig:Fundamental-Structures}, in particular, with
$$X_{u_1}=X_{u_8}, \ \ X_{u_2}=X_{u_7}, \ \ X_{u_3}=X_{u_6}, \ \ X_{u_4}=X_{u_5},$$ $$Y_{u_1}=Y_{u_2}=Y_{u_3}=Y_{u_4},$$
$$Y_{u_5}=Y_{u_6}=Y_{u_7}=Y_{u_8},$$
and
\begin{equation}\label{equation below-and-above}
|\partial(U)\cap \mathcal{A}_{\mathcal{D}_n}(u_i)|=|\partial(U)\cap \mathcal{B}_{\mathcal{D}_n}(u_i)|=\frac{n-3}{2}
\end{equation}
for $i\in [1,8]$, i.e.,
the number of edges of $\partial(U)$ which are above arcs with respect to $u_i$ and the number of edges of $\partial(U)$ which are below arcs with respect to $u_i$ are both $\frac{n-3}{2}$. Moreover,
no vertex lies in the interior of the cycle $C$,
where $C$ denotes any of the $4$-cycles $u_1u_2u_3u_4$,
$u_3u_4u_5u_6$, $u_5u_6u_7u_8$ and $u_7u_8u_1u_2$.
\begin{figure}[!ht]
\centering\includegraphics[scale=1.0]{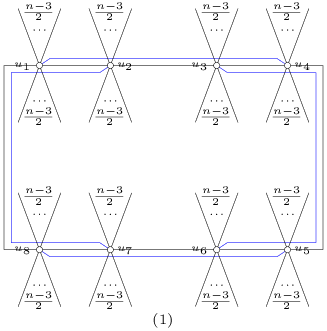} \hspace{60bp}
\centering\includegraphics[scale=1.0]{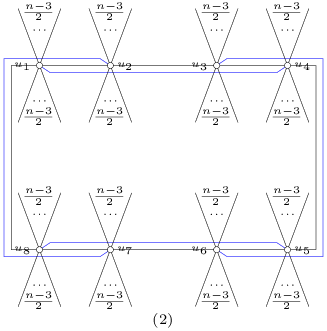} \caption{\small{
Fundamental structures}}\label{fig:Fundamental-Structures}
\end{figure}

Furthermore, if the drawing of the edge set $\mathcal{I}(U)$ in
$\mathcal{D}_n$ forms a fundamental structure, the edges of $E(U)$ are
called \verb"fundamental edges".
\end{definition}

\noindent $\bullet$ {\sl It can be verified  in Figure \ref{fig:gamma 5}  that there exist two fundamental structures in the drawing $\Gamma_5$, say $\mathcal{I}(U_1)$ and $\mathcal{I}(U_2)$, where $U_1$ and $U_2$ denote the set of eight vertices numbered as $2,10,26,18,22,30,14,6$ and numbered as $3,11,27,19,23,31,15,7$ respectively,

Now we are ready to characterize each step in
detail.

\noindent \textbf{Step 1}. Obtaining $\Gamma_n^*$ from $\Gamma_n$.

In this step, we construct the drawing $\Gamma_n^*$ satisfying the
following three properties:

\noindent {\rm Property 1:}   {\it For any vertex $u$ in $\Gamma_n^*$, the number of left arcs and the number of right arcs with respect to $u$ are almost equal,  precisely speaking, either
 \begin{equation}\label{equation left is more}
 (|\mathcal{L}_{\Gamma_n^*}(u)|,|\mathcal{R}_{\Gamma_n^*}(u)|)=(\frac{n+1}{2},\frac{n-1}{2})
  \end{equation}
 or
  \begin{equation}\label{equation left is less}
(|\mathcal{L}_{\Gamma_n^*}(u)|,|\mathcal{R}_{\Gamma_n^*}(u)|)=(\frac{n-1}{2},\frac{n+1}{2});
 \end{equation}}

\noindent {\rm Property 2:} {\it In the drawing $\Gamma_n^*$, there exist $2^{\frac{n-3}{2}}$ fundamental structures;}

\noindent {\rm Property 3:}  {\it For any $n\in \{5,7,9\}$, there exists a
decomposition of $V(Q_n)$ into several disjoint cycles which satisfy some conditions, (called enclosed cycles which will be introduced later)
say $\mathcal {C}_1,\mathcal {C}_2,\ldots,\mathcal {C}_m$, in
$\Gamma_n^*$ such that no fundamental edge is contained in any enclosed cycle, where $m\geq 1$. }

\medskip

The first thing worth noting is that the adjustments in this step change only the drawing of the adjusted edge within the `very small neighborhoods' of both its ends, in particular, the adjustments change only the nature of the adjusted edge with both its ends but keep the position relation between the adjusted edge and other edges unchanged from the
topological point of view.
To make it clear, we take $n=9$ for example in Figure \ref{fig:Adjustment-locally} to show the adjustments in this step. Suppose first in the drawing $\Gamma_9$, $|\mathcal{L}_{\mathcal{D}_n}(u)|=2$ and $|\mathcal{R}_{\mathcal{D}_n}(u)|=7$, i.e., the number of left arcs and the number of right arcs with respect $u$ are not almost equal. By redrawing the edges $e_3$ and $e_4$ within the very small neighborhood of $u$ to make them to be the left arcs with respect to $u$, we make the number of left arcs and the number of right arcs with respect to $u$ to be almost equal. As seen in Figure \ref{fig:Adjustment-locally}, the adjustments change only the natures of $e_3$ and $e_4$ with respect to $u$ (actually, the natures with respect to another end of $e_3,e_4$ also need to be changed in order to keep $e_3,e_4$ self-symmetric, here we emphasize the adjustment only with respect to $u$) but not change the topological position of $e_3,e_4$ with other edges.
 \begin{figure}[!ht]
\centering
\includegraphics[scale=1.0]{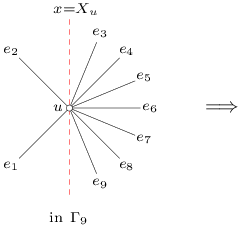}\hspace{15bp}
\includegraphics[scale=1.0]{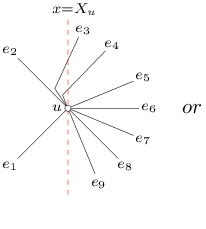}
\includegraphics[scale=1.0]{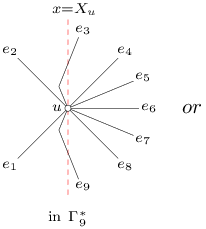}
\includegraphics[scale=1.0]{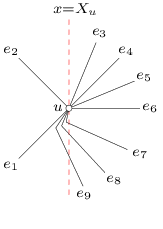}
\caption{\small{The characteristics of the adjustments in Step 1}}
\label{fig:Adjustment-locally}
\end{figure}
In fact, the selection of the edges to be adjusted is decided by elaborate attempts instead
of a general rule in the cases when $n\in \{5,7,9\}$.

For Property 2, as observed earlier, in the drawing $\Gamma_5$, there exist $2=2^{\frac{5-3}{2}}$ fundamental structures. To make $\Gamma_n^*$ satisfy Property 2, during the whole process of obtaining $\Gamma_{n+2}$ from $\Gamma_n$, we shall create two new fundamental structures  out of any fundamental structures in $\Gamma_n$, and moreover, during the process of obtaining $\Gamma_n^*$ from $\Gamma_n$, we shall preserve the fundamental structures, i.e., we should avoid adjusting any fundamental edge and changing the number of below arcs and the number of above arcs with respect to any vertex in some fundamental structure in $\Gamma_n^*$ (see \eqref{equation below-and-above}).

To illustrate Property 3, we need to give the definition of enclosed cycles.

\begin{definition}\label{definition enclosed-cycles} Let $\mathcal{C}$ be a cycle in $\Gamma_n^*$. We call $\mathcal{C}$ an {\bf enclosed cycle} provided that the  following two conditions hold.

\noindent {\rm Condition 1}:  Let $e_1,e_2$ be two adjacent edges in the cycle $\mathcal{C}$, say their common endpoint is $u$. Then one of the following types given in Figure \ref{fig:enclosedcycletwoedges} holds for the drawing within the very small neighborhood of $u$, where the bold lines represent the edges $e_1,e_2$.

\begin{figure}[!ht]
\centering
\includegraphics[scale=1.0]{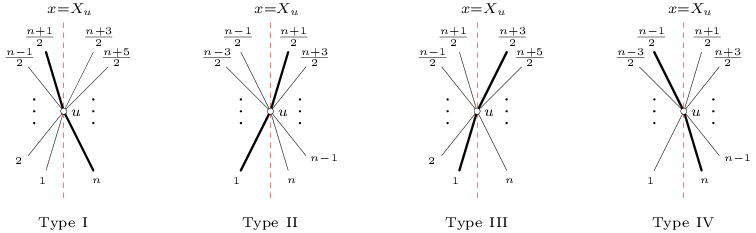}
\caption{\small{Two adjacent edges in the cycle $C$ incident with $u$}}
\label{fig:enclosedcycletwoedges}
\end{figure}

\noindent {\rm Condition 2}: Let $e=u_1u_2$ be an edge in the cycle $\mathcal{C}$.

Suppose $Y_{u_1}=Y_{u_2}$, say $X_{u_1}<X_{u_2}$. Then the relation of $e$ with
its both ends $u_1$ and $u_2$ is displayed as some diagram in Figure \ref{fig:edge-for-enclosed-cycle}.  In particular, $$|\mathcal{L}_{\Gamma_n^*}(u_1)|=|\mathcal{R}_{\Gamma_n^*}(u_2)|,$$ and $$e\in \mathcal{B}_{\Gamma_n^*}(u_1) \ \ \  (e\in \mathcal{A}_{\Gamma_n^*}(u_1))\ \ \mbox{ if and only if}\ \  e\in \mathcal{B}_{\Gamma_n^*}(u_2) \ \ \ (e\in \mathcal{A}_{\Gamma_n^*}(u_2)),$$
and $$e\in \mathcal{L}_{\Gamma_n^*}(u_1) \ \ \  (e\in \mathcal{R}_{\Gamma_n^*}(u_1))\ \ \mbox{ if and only if}\ \  e\in \mathcal{R}_{\Gamma_n^*}(u_2) \ \ \ (e\in \mathcal{L}_{\Gamma_n^*}(u_2)).$$
\begin{figure}[!ht]
\centering
\includegraphics[scale=0.64]{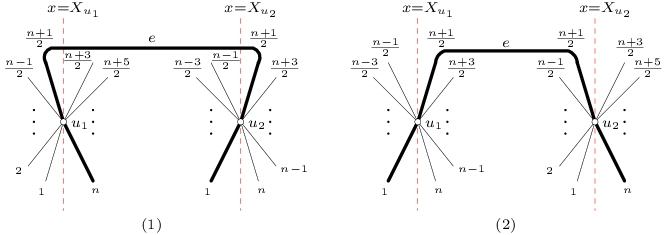}\hspace{10bp}
\includegraphics[scale=0.64]{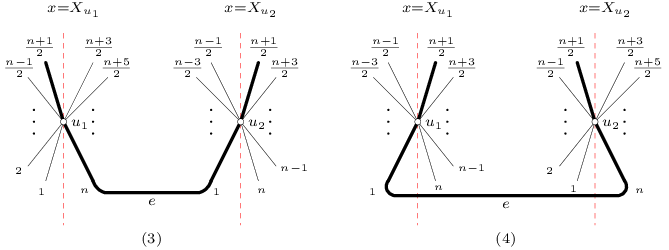}

\vspace{10bp}
\includegraphics[scale=0.64]{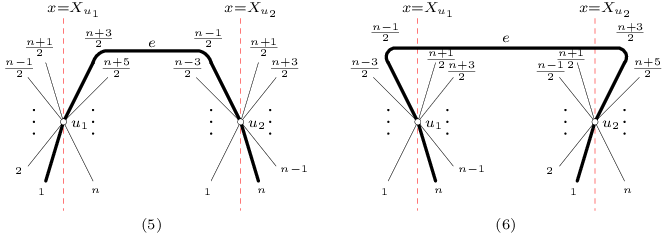}\hspace{10bp}
\includegraphics[scale=0.64]{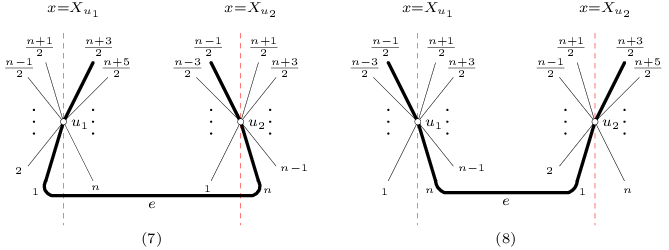}
\caption{\small{The edge $u_1u_2$ in an enclosed cycle with $Y_{u_1}=Y_{u_2}$}}\label{fig:edge-for-enclosed-cycle}
\end{figure}

Suppose $X_{u_1}=X_{u_2}$, say $Y_{u_1}>Y_{u_2}$. Then the relation of $e$ with
its both ends $u_1$ and $u_2$ is displayed as some diagram in Figure  \ref{fig:edge-for-enclosed-cycle2}. In particular, $$|\mathcal{L}_{\Gamma_n^*}(u_1)|=|\mathcal{L}_{\Gamma_n^*}(u_2)|,$$
and $$e\in \mathcal{B}_{\Gamma_n^*}(u_1) \ \ \  (e\in \mathcal{A}_{\Gamma_n^*}(u_1))\ \ \mbox{ if and only if}\ \  e\in \mathcal{A}_{\Gamma_n^*}(u_2) \ \ \ (e\in \mathcal{B}_{\Gamma_n^*}(u_2)),$$
and $$e\in \mathcal{L}_{\Gamma_n^*}(u_1) \ \ \  (e\in \mathcal{R}_{\Gamma_n^*}(u_1))\ \ \mbox{ if and only if}\ \  e\in \mathcal{L}_{\Gamma_n^*}(u_2) \ \ \ (e\in \mathcal{R}_{\Gamma_n^*}(u_2)).$$
\begin{figure}[!ht]
\centering
\includegraphics[scale=0.7]{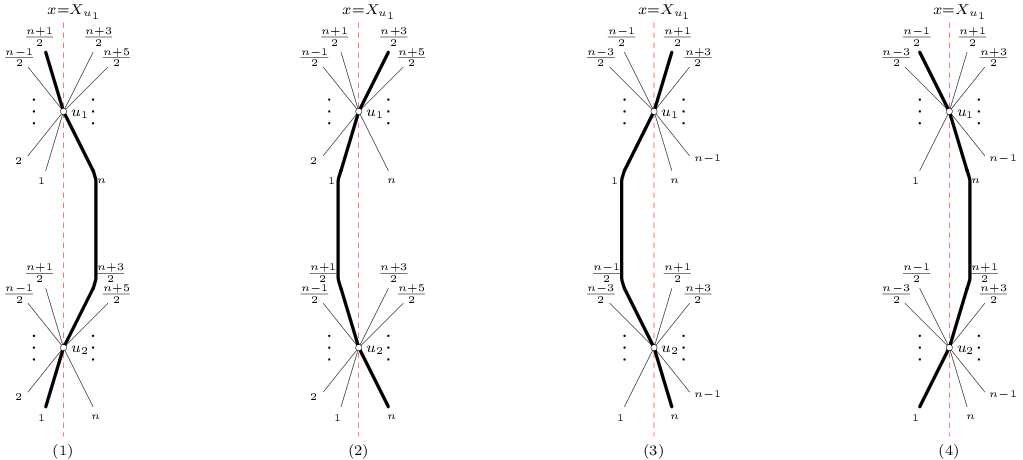}\vspace{15bp}
\includegraphics[scale=0.7]{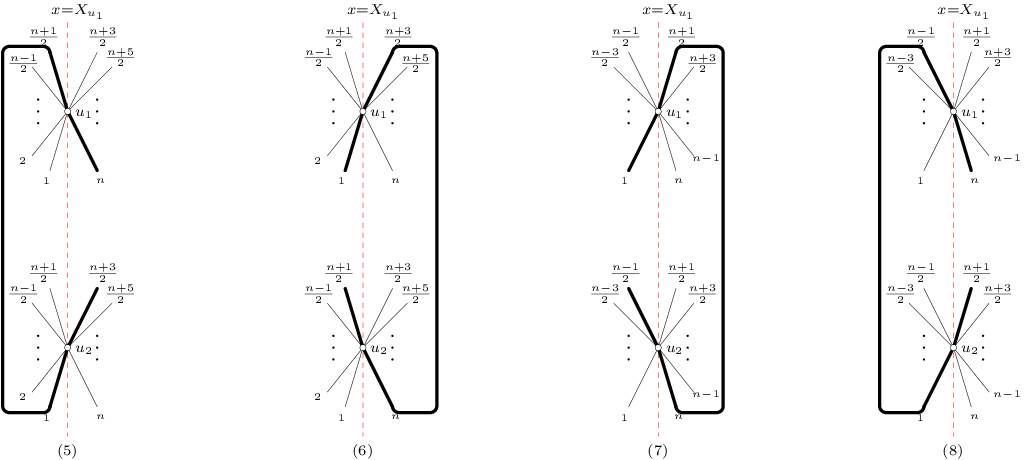}
\caption{\small{The edge $u_1u_2$ in an enclosed cycle with $X_{u_1}=X_{u_2}$}}\label{fig:edge-for-enclosed-cycle2}
\end{figure}
\end{definition}

Finding a decomposition of $V(Q_n)$ into several disjoint enclosed cycles is preparing, in fact, for decreasing the crossings within the neighborhoods of the small $4$-cycle $u^{(00)} u^{(10)}u^{(11)}u^{(01)}$ during the process from $\widetilde{\Gamma}_{n+2}$ to $\Gamma_{n+2}$, where $u$ denotes an arbitrary vertex in $\Gamma_n^*$ for $n\in \{5,7,9\}$, on which the specific operation will be displayed in Step 3. We can not find a general rule on how to construct $\Gamma_n^*$ from $\Gamma_n$ when $n\in \{5,7,9\}$. Instead, the process of constructing $\Gamma_5^*, \Gamma_7^*, \Gamma_9^*$ from $\Gamma_5, \Gamma_7, \Gamma_9$ is designed by elaborate attempts. Therefore, we shall display the process for $n\in \{5,7,9\}$ by figures. As for the case of $n\geq 11$, to avoid wasteful duplication in Step 2 and Step 3, we plan to leave the general rule on how to construct $\Gamma_n^*$ from $\Gamma_n$ to the later part of this subsection when we complete the description in Step 3.

From Figure \ref{fig:gamma 5} for the drawing $\Gamma_5$,
Figure \ref{fig:gamma 5star} for the drawing $\Gamma_5^*$,
Figure \ref{fig:wholedrawing7} (note that Figures \ref{fig:wholedrawing7}-- \ref{fig:Gamma9star} are placed at the end of the paper) for the drawing $\Gamma_7$, Figure \ref{fig:wholedrawing77} for the drawing $\Gamma_7^*$,
Figure \ref{fig:gamma 9} for the drawing $\Gamma_9$,
Figure \ref{fig:Gamma9star} for the drawing $\Gamma_9^*$,  one can see the adjustments of constructing $\Gamma_n^*$ from $\Gamma_n$ for $n\in \{5,7,9\}$. In Figure \ref{fig:gamma 5star}, Figure \ref{fig:wholedrawing77} and Figure \ref{fig:Gamma9star} for the drawings $\Gamma_5^*, \Gamma_7^*, \Gamma_9^*$, the edges in any enclosed cycle are marked in bold.
\begin{figure}[!ht]
\centering
\includegraphics[scale=1.0]{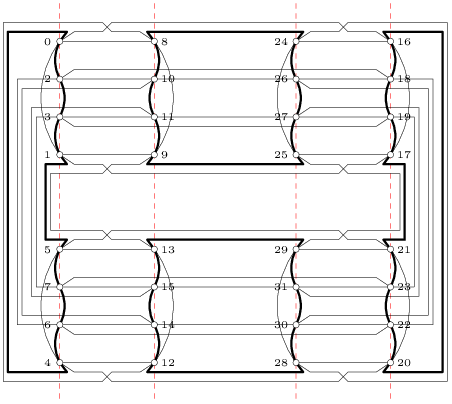}
\caption{\small{The drawing $\Gamma_5^*$}}\label{fig:gamma 5star}
\end{figure}

In Figure \ref{fig:gamma 9} for the drawing $\Gamma_9$, in
Figure \ref{fig:Gamma9star} for the drawing $\Gamma_9^*$, we depicted only a part of the whole drawing.  It is mainly because the page size is limited.  To make the adjustments to be seen, we display only the drawings of the edges incident
with the vertices $u$ with  $X_u\in\{-2,-1\}$ and $3\leq Y_u\leq 4$ in both Figures.
The whole drawing of $\Gamma_9^*$ ($\Gamma_9$) can be divided into eight identical parts (see Figure \ref{fig:eight-parts} for a schematic illustration), denoted $L_1,L_2,L_3,L_4,R_1,R_2,R_3,R_4$, among which, $L_1$ and $L_2$ are reflections of each other (including which edges are marked in bold) with the `mirror' $y=\frac{5}{2}$;  $L_3,L_4$ as a whole part is a reflection of $L_1,L_2$ as a whole part with the ` mirror' $y=0$; $R_1,R_2,R_3,R_4$ as a whole part is a reflection of $L_1,L_2,L_3,L_4$ as a whole part with the ` mirror' $x=0$. \begin{figure}[!ht]
\centering
\includegraphics[scale=0.8]{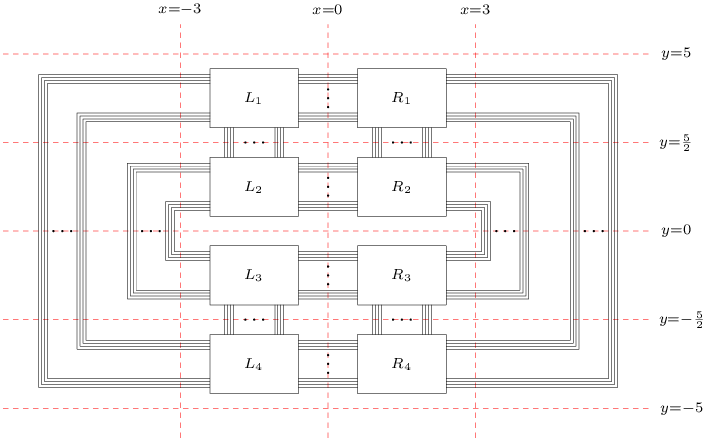}
\caption{\small{The relation of the eight identical parts composing the whole drawing}}\label{fig:eight-parts}
\end{figure}
If fact, each one of the drawings $\Gamma_5,\Gamma_5^*,\Gamma_7,\Gamma_7^*$ we give in Figures \ref{fig:gamma 5}, \ref{fig:gamma 5star}, \ref{fig:wholedrawing7}, \ref{fig:wholedrawing77} can be divided into eight identical parts in the way as stated above.
In a later part of this subsection, we shall illustrate our constructions in Figures \ref{fig:aid diagrams1} (1) and (2), \ref{fig:aid diagrams2},  \ref{fig:aid diagrams3},  \ref{fig:aid diagrams4}, \ref{fig:aid diagramsGamma6} and \ref{fig:Gamma8} by giving the corresponding drawings of $\Gamma_5, \Gamma_5^*, \widetilde{\Gamma}_7, \Gamma_7, \Gamma_7^*, \Gamma_6, \Gamma_8$ within the part $L_1$.

\noindent  \textbf{Step 2}. Obtaining $\widetilde{\Gamma}_{n+2}$
from $\Gamma_n^*$.

Let $u_1u_2$ be an arbitrary edge in the drawing $\Gamma_n^*$. By the inductive rule for the arrangement of vertices (see \eqref{equation
X for odd n}, \eqref{equation Y for odd n} and Figure \ref{fig:inductiveruleforvertices}), we can replace
$u_1u_2$  by a bunch of four new edges
$u_1^{(00)}u_2^{(00)},u_1^{(10)}u_2^{(10)},u_1^{(11)}u_2^{(11)},u_1^{(01)}u_2^{(01)}$
such that the bunch is drawn along
the original route of $u_1u_2$ in $\Gamma_n^*$. The local drawing (we call it `mesh-like structure') around the eight vertices $u_1^{(00)},u_1^{(10)},u_1^{(11)},u_1^{(01)}$  and $u_2^{(00)}, u_2^{(10)}, u_2^{(11)}, u_2^{(01)}$ in $\widetilde{\Gamma}_{n+2}$ are as follows. We take an edge $u_1^{(ab)}u_2^{(ab)}$ for example, where $ab\in \{00,10,11,01\}$.
In the drawing $\widetilde{\Gamma}_{n+2}$, the natures of the edge $u_1^{(ab)}u_2^{(ab)}$ with respect to its ends $u_1^{(ab)}$ and $u_2^{(ab)}$ are the same as the natures of the edge $u_1u_2$ with respect to vertices $u_1$ and $u_2$ in $\Gamma_n^*$.
More precisely,
$$u_1^{(ab)}u_2^{(ab)}\in \mathcal{L}_{\widetilde{\Gamma}_{n+2}}(u_i^{(ab)})$$ $$((u_1^{(ab)}u_2^{(ab)}\in \mathcal{R}_{\widetilde{\Gamma}_{n+2}}(u_i^{(ab)}), u_1^{(ab)}u_2^{(ab)}\in \mathcal{B}_{\widetilde{\Gamma}_{n+2}}(u_i^{(ab)}), u_1^{(ab)}u_2^{(ab)}\in \mathcal{A}_{\widetilde{\Gamma}_{n+2}}(u_i^{(ab)}), \mbox{ correspondingly})$$ provided that
$$u_1u_2\in \mathcal{L}_{\Gamma_n^*}(u_i)$$
$$((u_1u_2\in \mathcal{R}_{\Gamma_n^*}(u_i), u_1u_2\in \mathcal{B}_{\Gamma_n^*}(u_i), u_1u_2\in \mathcal{A}_{\Gamma_n^*}(u_i)).$$
Moreover, the edges $u_i^{(00)}u_i^{(10)}$,
$u_i^{(10)}u_i^{(11)}$, $u_i^{(11)}u_i^{(01)}$ are drawn precisely
along the line $x=X_{u_i}$. The edge $u_i^{(00)}u_i^{(01)}$ is
drawn to be an arc on the right side or on the left side of the line
$x=X_{u_i}$ according to  \eqref{equation left is more} or \eqref{equation left is less}
holds for $u_i$ in $\Gamma_n^*$ respectively,
where $i=1,2$.

In Figure \ref{fig:local drawing}, Diagram (1) and Diagram (2) depict the mesh-like structure formed by edges incident with $u^{(00)},u^{(10)},u^{(11)}, u^{(01)}$ within the neighborhood of the small $4$-cycle $u^{(00)}u^{(10)}u^{(11)}u^{(01)}$ in $\widetilde{\Gamma}_{n+2}$ for the case of $(|\mathcal{L}_{\Gamma_n^*}(u)|,|\mathcal{R}_{\Gamma_n^*}(u)|)=(\frac{n+1}{2},\frac{n-1}{2})$
and for the case of $(|\mathcal{L}_{\Gamma_n^*}(u)|,|\mathcal{R}_{\Gamma_n^*}(u)|)=(\frac{n-1}{2},\frac{n+1}{2})$
respectively, where $u$ denotes an arbitrary vertex in $\Gamma_n^*$.
\begin{figure}[!ht]
\centering
\includegraphics[scale=1.0]{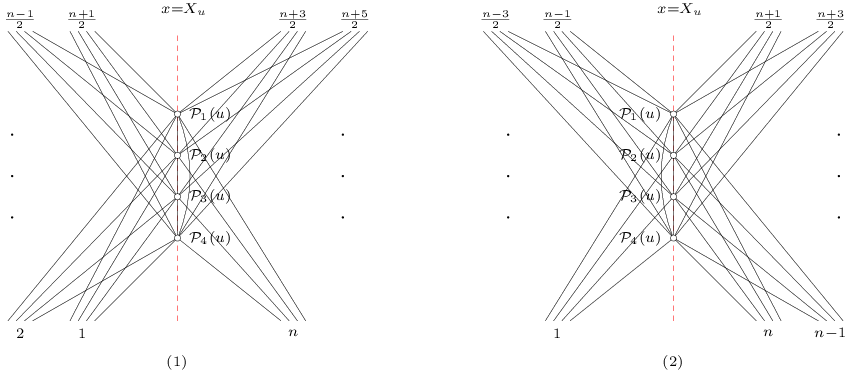}
\caption{\small{The mesh-like structure around $u^{(00)}u^{(10)}u^{(11)}u^{(01)}$ in $\widetilde{\Gamma}_{n+2}$ constructed from  $\Gamma_n^*$}}\label{fig:local drawing}
\end{figure}

On the other hand, during the process in Step 2 we see that the drawing $\widetilde{\Gamma}_{n+2}$ would inherit the characteristic of every edge being self-symmetric from $\Gamma_n^*$. Moreover, the four edges
$u_1^{(00)}u_2^{(00)},
u_1^{(10)}u_2^{(10)},
u_1^{(11)}u_2^{(11)},
u_1^{(01)}u_2^{(01)}$
extended from $u_1u_2$ {\bf do not cross each other}.
This can be seen from Figure \ref{fig:keep-selfsymmetric}, in which we give the corresponding drawings of the four edges
$u_1^{(00)}u_2^{(00)},
u_1^{(10)}u_2^{(10)},
u_1^{(11)}u_2^{(11)},
u_1^{(01)}u_2^{(01)}$ by taking the cases that the self-symmetric edge $u_1u_2$ is depicted as
Diagrams (1), (4) in Figure \ref{fig:Selfsymmx=x} and Diagrams (1), (3) in Figure \ref{fig:Selfsymmy=y} when $Y_{\widehat{u_1}}<Y_{u_1}$ for examples.
\begin{figure}[!ht]
\centering
\includegraphics[scale=0.75]{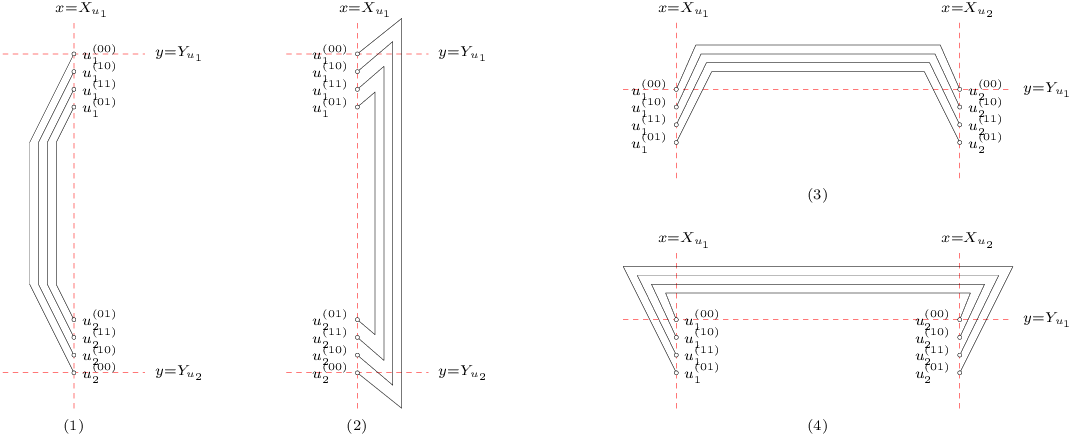}
\caption{\small{$\widetilde{\Gamma}_{n+2}$ inherits the characteristic of every edge being self-symmetric from $\Gamma_n^*$}}\label{fig:keep-selfsymmetric}
\end{figure}

The process of obtaining $\widetilde{\Gamma}_7$ from $\Gamma_5^*$
shown in Figures \ref{fig:aid diagrams1} and \ref{fig:aid diagrams2} will be helpful for us to
understand accurately the adjustments described in this step.

\noindent \textbf{Step 3}. Obtaining $\Gamma_{n+2}$ from
$\widetilde{\Gamma}_{n+2}$.

To obtain  the drawing $\Gamma_{n+2}$, we need to make two kinds of
adjustments on the edges in $\widetilde{\Gamma}_{n+2}$. The first
kind is applied only on the edges in  $\widetilde{\Gamma}_7,\widetilde{\Gamma}_9,\widetilde{\Gamma}_{11}$, which are associated with the
enclosed cycles in $\Gamma_5^*,\Gamma_7^*,\Gamma_9^*$, respectively. The second is applied on the edges in $\widetilde{\Gamma}_{n+2}$ which are associated with the fundamental structures in
$\Gamma_n^*$ for all $n\geq 5$.

Now we take the first kind of adjustments and suppose $n\in \{5,7,9\}$ at first. Let $u$ be an arbitrary vertex in the drawing $\Gamma_n^*$. Recall that the drawing $\Gamma_n^*$ has Property 3, i.e., $u$ is in some enclosed cycle $\mathcal{C}$ (see Figure \ref{fig:enclosedcycletwoedges} for the four types of the local drawing within the very small neighborhood of $u$ in $\Gamma_n^*$).  Combined with  the adjustment given in Step 2 (see Figure \ref{fig:local drawing}), we have that the mesh-like structures around the small $4$-cycle $u^{(00)}u^{(10)}u^{(11)}u^{(01)}$ are depicted as in Figure \ref{fig:basic-mesh-bold} of which Diagrams (1)-(4) correspond to Types {\rm I}-{\rm IV} of Figure \ref{fig:enclosedcycletwoedges}, where
the two bunches marked in bold are the group of edges extended from the edges of some enclosed cycle in $\Gamma_n^*$.
\begin{figure}[!ht]
\centering
\includegraphics[scale=0.75]{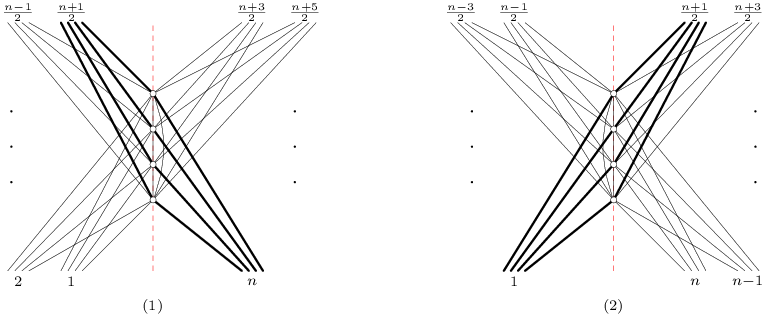}

\vspace{15bp}
\includegraphics[scale=0.75]{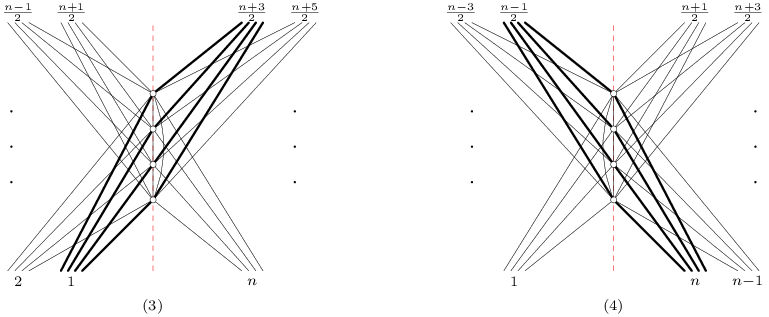}
\caption{\small{Mesh-like structures with depicting the edges extended from enclosed cycles}}\label{fig:basic-mesh-bold}
\end{figure}

Corresponding to each mesh-like structure depicted in Figure \ref{fig:basic-mesh-bold}, we displayed in Figure \ref{fig:basic-mesh-bold-adjusted} the {\bf adjusted mesh-like structure} around $u^{(00)}u^{(10)}u^{(11)}u^{(01)}$ after the first kind of adjustments done. In particular, we reversed two edges from those two bunches within the neighborhoods of the small $4$-cycle $u^{(00)}u^{(10)}u^{(11)}u^{(01)}$ and marked the two reversed edges in bold.
\begin{figure}[!ht]
\centering
\includegraphics[scale=0.75]{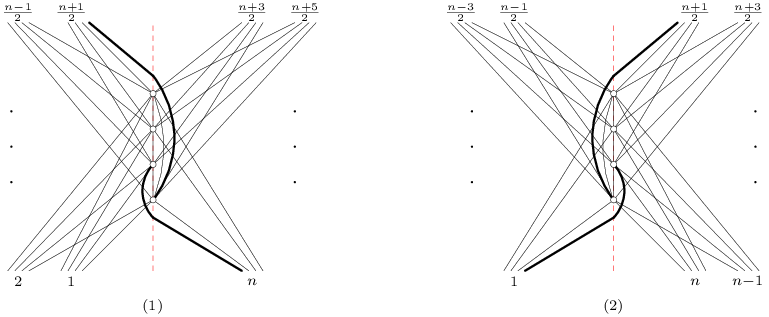}

\vspace{15bp}
\includegraphics[scale=0.75]{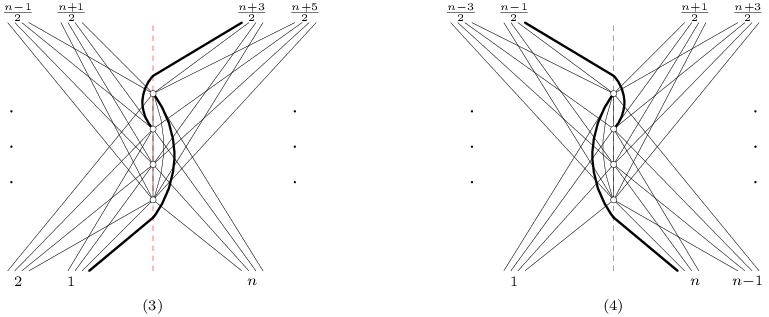}
\caption{\small{Adjusted mesh-like structures around $u^{(00)}u^{(10)}u^{(11)}u^{(0)}$}}
\label{fig:basic-mesh-bold-adjusted}
\end{figure}

Combined with Condition 2 (see Figure \ref{fig:edge-for-enclosed-cycle} and Figure \ref{fig:edge-for-enclosed-cycle2}) in the definition of enclosed cycles, Definition \ref{definition enclosed-cycles}, we can derive that the adjusted edges still keep self-symmetric, and that the four edges in any bunch which is adjusted, {\bf do not cross each other} and still are in a bunch along their original route before the adjustments done, what differ is the relative position of the four edges changed `a little'. We show this by Figure \ref{fig:Extended-enclosed-cycle1} and Figure \ref{fig:Extended-enclosed-cycle2}, in which each diagram corresponds to each diagram in  Figure \ref{fig:edge-for-enclosed-cycle} and Figure \ref{fig:edge-for-enclosed-cycle2} respectively, where in Figure \ref{fig:Extended-enclosed-cycle1} and Figure \ref{fig:Extended-enclosed-cycle2} the adjusted edges are marked in bold.
\begin{figure}[!ht]
\centering
\includegraphics[scale=0.75]{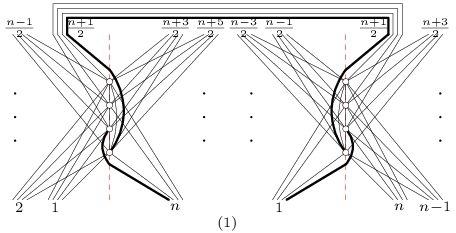}\hspace{50bp}
\includegraphics[scale=0.75]{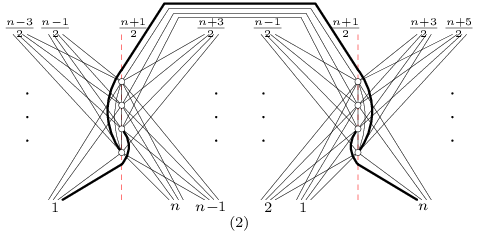}

\vspace{15bp}
\includegraphics[scale=0.75]{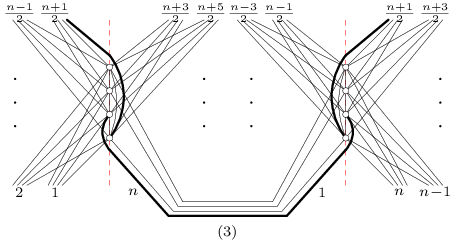}\hspace{50bp}
\includegraphics[scale=0.75]{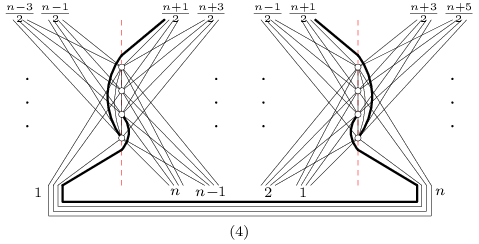}

\vspace{15bp}
\includegraphics[scale=0.75]{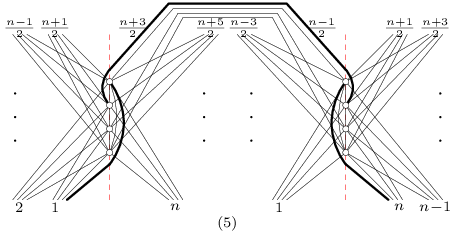}\hspace{50bp}
\includegraphics[scale=0.75]{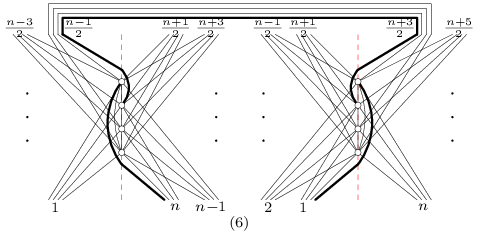}

\vspace{15bp}
\includegraphics[scale=0.75]{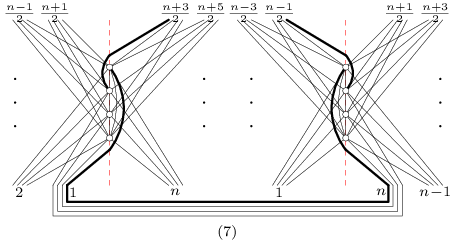}\hspace{50bp}
\includegraphics[scale=0.75]{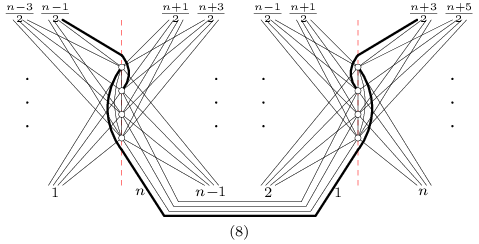}
\caption{\small{The adjustment with edges in enclosed cycles I }}\label{fig:Extended-enclosed-cycle1}
\end{figure}

\begin{figure}[!ht]
\centering
\includegraphics[scale=0.75]{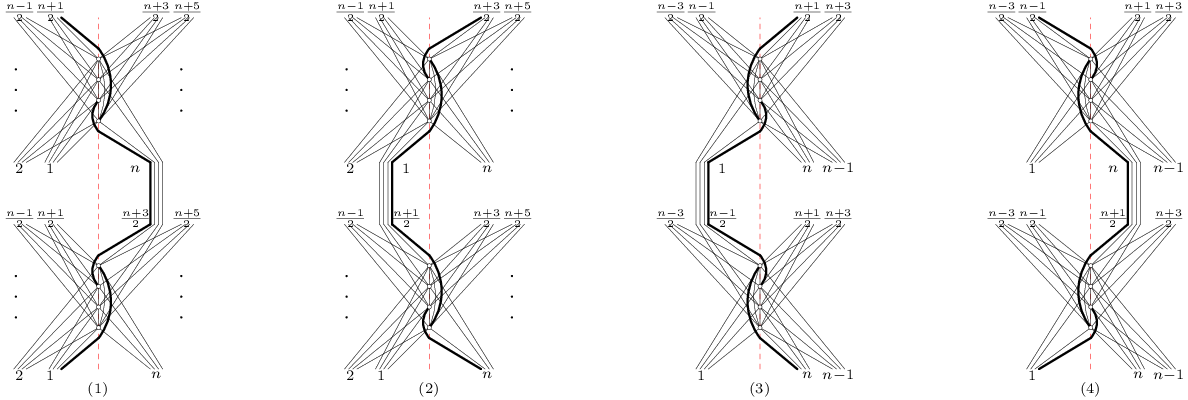}

\vspace{15bp}
\includegraphics[scale=0.75]{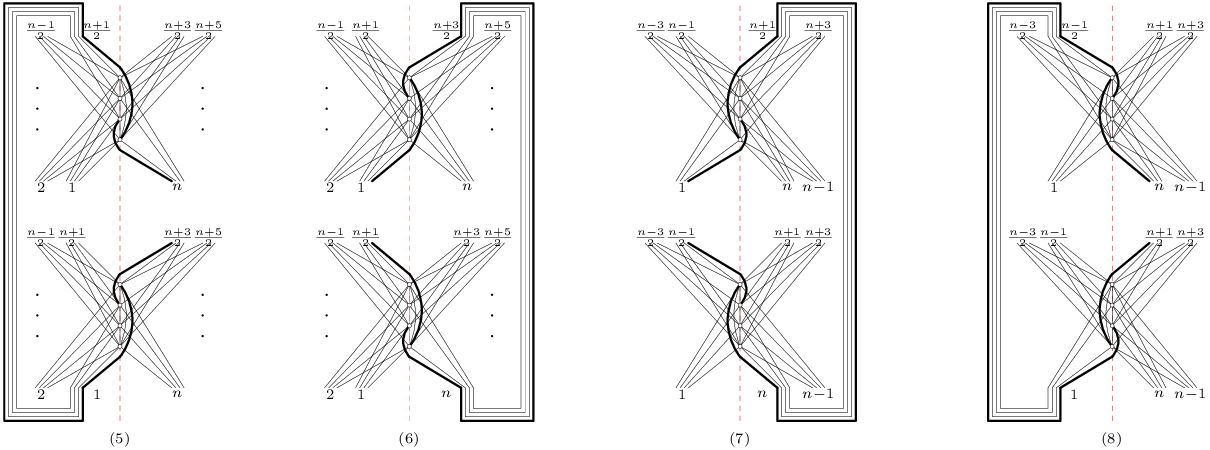}
\caption{\small{The adjustment with edges in enclosed cycles II }}\label{fig:Extended-enclosed-cycle2}
\end{figure}

Next we take the second kind of adjustments in this step. This kind of adjustments are applied on edges in $\widetilde{\Gamma}_{n+2}$ for all odd integers $n\geq 5$, which are associated with the fundamental structures in $\Gamma_n^*$ and will change the route of the adjusted edge `globally'. Take a
vertex subset $U=\{u_1,u_2,\ldots,u_8\}$ of $V(Q_n)$ such that the
drawing of $\mathcal{I}(U)$ forms a fundamental structure in
$\Gamma_n^*$, which is drawn without loss of generality as Diagram (1) of
Figure \ref{fig:Fundamental-Structures}.
\begin{figure}[!ht]
\centering\includegraphics[scale=0.9]{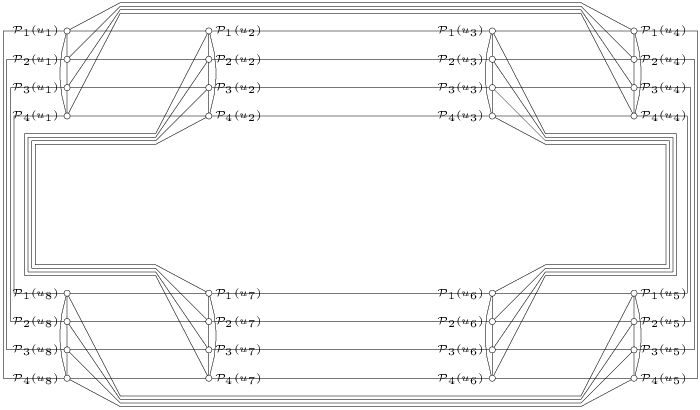} \caption{\small{The
drawing of $E(\pi(U))$ in
$\widetilde{\Gamma}_{n+2}$}}\label{fig:BeforeFundamental}
\end{figure}
We see that the edges of
$E(\pi(U))$ in $\widetilde{\Gamma}_{n+2}$ will be drawn as in Figure
\ref{fig:BeforeFundamental}, where
$$\pi(U)=\bigcup\limits_{i=1}^8\{u_i^{(00)},u_i^{(10)},u_i^{(11)},u_i^{(01)}\}.$$
Recall that no fundamental edge belongs to any enclosed cycle in $\Gamma_5^*,\Gamma_7^*,\Gamma_9^*$ given in Property 3. On the other hand, the first kind of adjustments changed only the edges from the bunches extended from the edges of enclosed cycles in $\Gamma_n^*$ with $n\in \{5,7,9\}$. Hence, the drawing for the edges of
$E(\pi(U))$ after the first kinds of adjustment on $\widetilde{\Gamma}_{n+2}$ keeps the same as in Figure \ref{fig:BeforeFundamental} for $n\in\{5,7,9\}$. Then we adjust the following
edges
\begin{eqnarray}
\mathcal{P}_2(u_1)\mathcal{P}_2(u_4),
\mathcal{P}_3(u_1)\mathcal{P}_3(u_4),
\mathcal{P}_4(u_1)\mathcal{P}_4(u_4),\nonumber \\
\mathcal{P}_1(u_3)\mathcal{P}_4(u_6),
\mathcal{P}_2(u_3)\mathcal{P}_3(u_6),
\mathcal{P}_3(u_3)\mathcal{P}_2(u_6),\nonumber \\
\mathcal{P}_1(u_5)\mathcal{P}_1(u_8),
\mathcal{P}_2(u_5)\mathcal{P}_2(u_8),
\mathcal{P}_3(u_5)\mathcal{P}_3(u_8),\nonumber \\
\mathcal{P}_1(u_2)\mathcal{P}_4(u_7),
\mathcal{P}_2(u_2)\mathcal{P}_3(u_7),
\mathcal{P}_3(u_2)\mathcal{P}_2(u_7),\nonumber
\end{eqnarray}
which are shown as in Figure \ref{fig:FinalFudametal}.  It is worthwhile to note during the process of the second kind adjustments, two new fundamental structures were created from the old one, which are emphasized in blue. Moreover, by \eqref{equation below-and-above}, we can verify that
$$\nu_{\Gamma_{n+2}}\big{(}E(\pi(U)),\partial(\pi(U))\big{)}=\nu_{\widetilde{\Gamma}_{n+2}}\big{(}E(\pi(U)),\partial(\pi(U))\big{)},$$
and thus,
\begin{equation}\label{equation crossings inner and outer}
\nu_{\Gamma_{n+2}}\big{(}E(\pi(U)),\ E(Q_{n+2})\setminus E(\pi(U))
\big{)}=\nu_{\widetilde{\Gamma}_{n+2}}\big{(}E(\pi(U)),\
E(Q_{n+2})\setminus E(\pi(U))\big{)},
\end{equation}
meanwhile, we can verify from Figure \ref{fig:BeforeFundamental} and Figure \ref{fig:FinalFudametal} that
\begin{equation}\label{equation crossings inner crossings same}
\nu_{\Gamma_{n+2}}\big{(}E(\pi(U))\big{)}=\nu_{\widetilde{\Gamma}_{n+2}}\big{(}E(\pi(U))\big{)}-8.
\end{equation}

\begin{figure}[!ht]
\centering\includegraphics[scale=1.0]{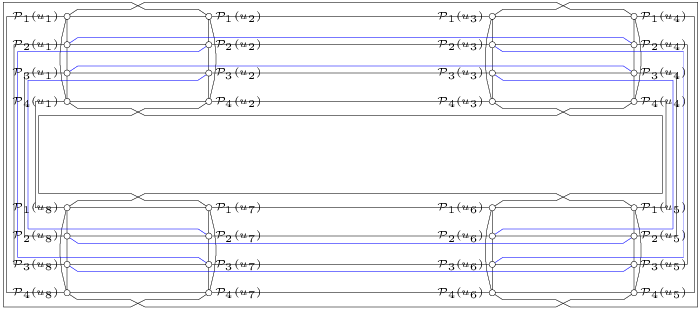} \caption{\small{The
drawing of $E(\pi(U))$ in $\Gamma_{n+2}$ depicting the adjustments associated with fundamental-structures}}\label{fig:FinalFudametal}
\end{figure}
This completes the description of the process in Step 3. To help the reader to understand the process described in the above three steps, we give Figures \ref{fig:aid diagrams1}, \ref{fig:aid diagrams2}, \ref{fig:aid diagrams3} and \ref{fig:aid diagrams4} to depict a complete process from $\Gamma_5$ to $\Gamma_7^*$.

\begin{figure}
\centering
\includegraphics[scale=1.3]{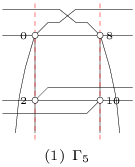} \hspace{50pt}
\includegraphics[scale=1.3]{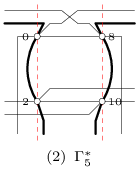}
\caption{\small{Auxiliary drawings illustrating the process from
$\Gamma_5$ to $\Gamma_5^*$}}\label{fig:aid diagrams1}
\end{figure}

\begin{figure}
\centering
\includegraphics[scale=1.15]{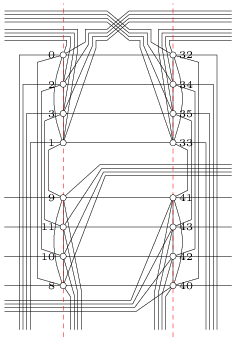}
\caption{\small{A part of $\widetilde{\Gamma}_7$} illustrating the process from
$\Gamma_5^*$ to $\widetilde{\Gamma}_7$}\label{fig:aid diagrams2}
\end{figure}

\begin{figure}
\centering
\includegraphics[scale=1.15]{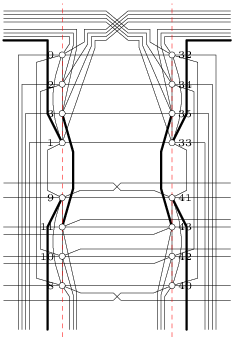}
\caption{\small{A part of $\Gamma_7$ illustrating the process from $\widetilde{\Gamma}_7$ to $\Gamma_7$ }}\label{fig:aid diagrams3}
\end{figure}

\begin{figure}
\centering
\includegraphics[scale=1.15]{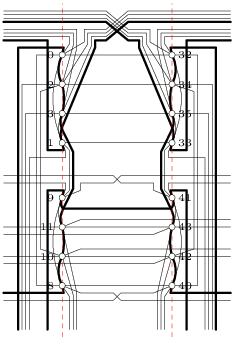}
\caption{\small{A part of $\Gamma_7^*$} illustrating the process from
$\Gamma_7$ to $\Gamma_7^*$}\label{fig:aid diagrams4}
\end{figure}

To complete the whole process, as stated in Step 1,  we still need to describe the general rule on how to obtain the desired drawing $\Gamma_n^*$ from $\Gamma_n$ when $n\geq 11$. Yet, there exist some differences between the process from $\Gamma_{11}$ to $\Gamma_{11}^*$ and the process from $\Gamma_n$ to $\Gamma_n^*$ for $n\geq 13$.

We first describe the rule for obtaining $\Gamma_{11}^*$ from $\Gamma_{11}$. To do this, we  need to find a perfect matching, denoted $\mathcal{M}$, in $\Gamma_9^*$, which meets the following:

\noindent {\bf Requirement for the matching $\mathcal{M}$}:  ``Let $u_1u_2$ be an arbitrary edge of the matching $\mathcal{M}$. Then $u_1u_2$ is the edge of some enclosed cycle with $$(|\mathcal{L}_{\Gamma_9^*}(u_i)|,|\mathcal{R}_{\Gamma_9^*}(u_i)|)=(4,5)\ \ \ \mbox{ or }\ \ \ (|\mathcal{L}_{\Gamma_9^*}(u_i)|,|\mathcal{R}_{\Gamma_9^*}(u_i)|)=(5,4)$$ according to $$u_1u_2\in \mathcal{L}_{\Gamma_9^*}(u_i)\ \ \ \mbox{ or } \ \ \  u_1u_2\in \mathcal{R}_{\Gamma_9^*}(u_i)$$ respectively, where $i=1,2$."

Since every vertex belongs to some enclosed cycle in $\Gamma_9^*$,  by the definition of enclosed cycle (see Figure \ref{fig:edge-for-enclosed-cycle} and Figure \ref{fig:edge-for-enclosed-cycle2}), we derive that there exists a unique such perfect matching $\mathcal{M}$ in $\Gamma_9^*$. Let $u_1u_2$ be an arbitrary edge of the matching $\mathcal{M}$. We also see that
the relative position of the edge $u_1u_2$ with its both ends $u_1$ and $u_2$ is depicted as Diagrams (3)-(6) in Figure \ref{fig:edge-for-enclosed-cycle} or Diagrams (1), (3), (6), (8) in Figure \ref{fig:edge-for-enclosed-cycle2}. Hence, after the first kind of adjustments on the drawing $\widetilde{\Gamma}_{11}$, the obtained drawings for $u_1^{(00)}u_2^{(00)},u_1^{(10)}u_2^{(10)},u_1^{(11)}u_2^{(11)},u_1^{(01)}u_2^{(01)}$ are depicted as the corresponding Diagrams (3)-(6) in Figure \ref{fig:Extended-enclosed-cycle1} or Diagrams (1), (3), (6), (8) in Figure \ref{fig:Extended-enclosed-cycle2}. By Property 3, no fundamental edge belongs to any enclosed cycle in $\Gamma_9^*$, we have that $u_1u_2$ is not a fundamental edge. Although there is a possibility that $u_i$ is an end of some fundamental edge in $\Gamma_9^*$,
the second kind of adjustments on $\widetilde{\Gamma}_{11}$ do not change the nature whether an edge $e\in \mathcal{I}(u_i^{(ab)})$
is a left arc or a right arc with respect to $u_i^{(ab)}$, and naturally do not change the number of left arcs and the number of right arcs with respect to $u_i^{(ab)}$,
(see Figure \ref{fig:FinalFudametal}), where $i=1,2$ and $ab\in \{00,10,11,01\}$. Moreover, observe that the second kind of adjustments do not involve any one of the following edges  $$u_1^{(00)}u_1^{(10)},u_1^{(10)}u_1^{(11)},u_1^{(11)}u_1^{(01)},$$
$$u_2^{(00)}u_2^{(10)},u_2^{(10)}u_2^{(11)},u_2^{(11)}u_2^{(01)},$$
$$u_1^{(00)}u_2^{(00)},u_1^{(10)}u_2^{(10)},u_1^{(11)}u_2^{(11)},u_1^{(01)}u_2^{(01)}.$$

Now by taking the cases depicted as Diagram (3) in Figure \ref{fig:Extended-enclosed-cycle1} or as Diagram (1) in Figure \ref{fig:Extended-enclosed-cycle2} when the first kind of adjustments are taken on $\widetilde{\Gamma}_{11}$ for example, we show the rule on how to construct $\Gamma_{11}^*$ from $\Gamma_{11}$. As stated above,  although in both cases, the number of left arcs and the number of right arcs with respect to $u_i^{(ab)}$ in $\Gamma_{11}$ is the same as given by Diagram (3) in Figure \ref{fig:Extended-enclosed-cycle1} or as given by Diagram (1) in Figure \ref{fig:Extended-enclosed-cycle2} (the drawings before the second kind of adjustments are taken), the mesh-like structure around $u_i^{(00)}u_i^{(10)}u_i^{(11)}u_i^{(01)}$ has probably changed. We shall also see that the possible changes on the mesh-like structures do not affect our following effort to balance the number of left arcs and the number of right arcs. Hence, we shall illustrate only some edges of which the positions are definitely unaffected by the second kind of adjustments rather than giving the whole mesh-like structures.

Consider the case depicted as Diagram (3) in Figure \ref{fig:Extended-enclosed-cycle1}. We check first that in $\Gamma_{11}$, $$(|\mathcal{L}_{{\Gamma}_{11}}(\mathcal{P}_1(u_1))|,|\mathcal{R}_{{\Gamma}_{11}}(\mathcal{P}_1(u_1))|)=(\frac{9+1}{2},\frac{9+1}{2})=(5,5),$$
$$(|\mathcal{L}_{{\Gamma}_{11}}(\mathcal{P}_2(u_1))|,|\mathcal{R}_{{\Gamma}_{11}}(\mathcal{P}_2(u_1))|)=(\frac{9+1}{2},\frac{9-1}{2})=(5,4),$$
$$(|\mathcal{L}_{{\Gamma}_{11}}(\mathcal{P}_3(u_1))|,|\mathcal{R}_{{\Gamma}_{11}}(\mathcal{P}_3(u_1))|)=(\frac{9+3}{2},\frac{9-3}{2})=(6,3),$$
$$(|\mathcal{L}_{{\Gamma}_{11}}(\mathcal{P}_4(u_1))|,|\mathcal{R}_{{\Gamma}_{11}}(\mathcal{P}_4(u_1))|)=(\frac{9-1}{2},\frac{9+3}{2})=(4,6),$$
and that
$$(|\mathcal{L}_{{\Gamma}_{11}}(\mathcal{P}_1(u_2))|,|\mathcal{R}_{{\Gamma}_{11}}(\mathcal{P}_1(u_2))|)=(\frac{9+1}{2},\frac{9+1}{2})=(5,5),$$
$$(|\mathcal{L}_{{\Gamma}_{11}}(\mathcal{P}_2(u_2))|,|\mathcal{R}_{{\Gamma}_{11}}(\mathcal{P}_2(u_2))|)=(\frac{9-1}{2},\frac{9+1}{2})=(4,5),$$
$$(|\mathcal{L}_{{\Gamma}_{11}}(\mathcal{P}_3(u_2))|,|\mathcal{R}_{{\Gamma}_{11}}(\mathcal{P}_3(u_2))|)=(\frac{9-3}{2},\frac{9+3}{2})=(3,6),$$
$$(|\mathcal{L}_{{\Gamma}_{11}}(\mathcal{P}_4(u_2))|,|\mathcal{R}_{{\Gamma}_{11}}(\mathcal{P}_4(u_2))|)=(\frac{9+3}{2},\frac{9-1}{2})=(6,4).$$
Then we make the number of left arcs and the number of right arcs  with respect to  each one of $\mathcal{P}_1(u_1)$, $\mathcal{P}_2(u_1)$, $\mathcal{P}_3(u_1)$, $\mathcal{P}_4(u_1)$,
$\mathcal{P}_1(u_2)$, $\mathcal{P}_2(u_2)$,
$\mathcal{P}_3(u_2)$, $\mathcal{P}_4(u_2)$
almost equal
by distorting seven edges (see Figure \ref{fig:gamma11togamma*h}), $$\mathcal{P}_1(u_1)\mathcal{P}_2(u_1), \mathcal{P}_2(u_1)\mathcal{P}_3(u_1), \mathcal{P}_3(u_1)\mathcal{P}_4(u_1),$$
$$\mathcal{P}_1(u_2)\mathcal{P}_2(u_2), \mathcal{P}_2(u_2)\mathcal{P}_3(u_2), \mathcal{P}_3(u_2)\mathcal{P}_4(u_2),$$
$$\mathcal{P}_4(u_1)\mathcal{P}_4(u_2),$$
`locally', and check that in ${\Gamma}_{11}^*$, $$(|\mathcal{L}_{{\Gamma}_{11}^*}(\mathcal{P}_1(u_1))|,|\mathcal{R}_{{\Gamma}_{11}^*}(\mathcal{P}_1(u_1))|)=(\frac{9+1}{2},\frac{9+3}{2})=(5,6),$$
$$(|\mathcal{L}_{{\Gamma}_{11}^*}(\mathcal{P}_2(u_1))|,|\mathcal{R}_{{\Gamma}_{11}^*}(\mathcal{P}_2(u_1))|)=(\frac{9+1}{2},\frac{9+3}{2})=(5,6),$$
$$(|\mathcal{L}_{{\Gamma}_{11}^*}(\mathcal{P}_3(u_1))|,|\mathcal{R}_{{\Gamma}_{11}^*}(\mathcal{P}_3(u_1))|)=(\frac{9+3}{2},\frac{9+1}{2})=(6,5),$$
$$(|\mathcal{L}_{{\Gamma}_{11}^*}(\mathcal{P}_4(u_1))|,|\mathcal{R}_{{\Gamma}_{11}^*}(\mathcal{P}_4(u_1))|)=(\frac{9+1}{2},\frac{9+3}{2})=(5,6),$$
and
$$(|\mathcal{L}_{{\Gamma}_{11}^*}(\mathcal{P}_1(u_2))|,|\mathcal{R}_{{\Gamma}_{11}^*}(\mathcal{P}_1(u_2))|)=(\frac{9+3}{2},\frac{9+1}{2})=(6,5),$$
$$(|\mathcal{L}_{{\Gamma}_{11}^*}(\mathcal{P}_2(u_2))|,|\mathcal{R}_{{\Gamma}_{11}^*}(\mathcal{P}_2(u_2))|)=(\frac{9+3}{2},\frac{9+1}{2})=(6,5),$$
$$(|\mathcal{L}_{{\Gamma}_{11}^*}(\mathcal{P}_3(u_2))|,|\mathcal{R}_{{\Gamma}_{11}^*}(\mathcal{P}_3(u_2))|)=(\frac{9+1}{2},\frac{9+3}{2})=(5,6),$$
$$(|\mathcal{L}_{{\Gamma}_{11}^*}(\mathcal{P}_4(u_2))|,|\mathcal{R}_{{\Gamma}_{11}^*}(\mathcal{P}_4(u_2))|)=(\frac{9+3}{2},\frac{9+1}{2})=(6,5).$$
\begin{figure}[!ht]
\centering\includegraphics[scale=1.0]{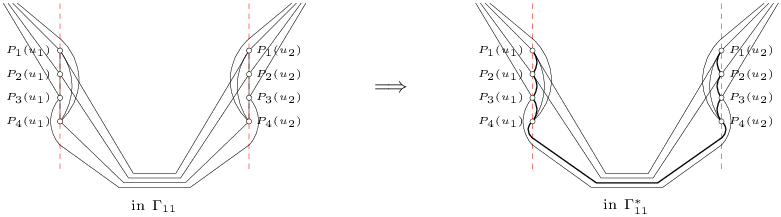} \caption{\small{Schematic diagram for the adjustments from $\Gamma_{11}$ to $\Gamma_{11}^*$ }}\label{fig:gamma11togamma*h}
\end{figure}

Consider the case depicted as Diagram (1) in Figure \ref{fig:Extended-enclosed-cycle2}. We check first that in $\Gamma_{11}$, $$(|\mathcal{L}_{{\Gamma}_{11}}(\mathcal{P}_1(u_1))|,|\mathcal{R}_{{\Gamma}_{11}}(\mathcal{P}_1(u_1))|)=(\frac{9+1}{2},\frac{9+1}{2})=(5,5),$$
$$(|\mathcal{L}_{{\Gamma}_{11}}(\mathcal{P}_2(u_1))|,|\mathcal{R}_{{\Gamma}_{11}}(\mathcal{P}_2(u_1))|)=(\frac{9+1}{2},\frac{9-1}{2})=(5,4),$$
$$(|\mathcal{L}_{{\Gamma}_{11}}(\mathcal{P}_3(u_1))|,|\mathcal{R}_{{\Gamma}_{11}}(\mathcal{P}_3(u_1))|)=(\frac{9+3}{2},\frac{9-3}{2})=(6,3),$$
$$(|\mathcal{L}_{{\Gamma}_{11}}(\mathcal{P}_4(u_1))|,|\mathcal{R}_{{\Gamma}_{11}}(\mathcal{P}_4(u_1))|)=(\frac{9-1}{2},\frac{9+3}{2})=(4,6),$$
and that
$$(|\mathcal{L}_{{\Gamma}_{11}}(\mathcal{P}_1(u_2))|,|\mathcal{R}_{{\Gamma}_{11}}(\mathcal{P}_1(u_2))|)=(\frac{9-1}{2},\frac{9+3}{2})=(4,6),$$
$$(|\mathcal{L}_{{\Gamma}_{11}}(\mathcal{P}_2(u_2))|,|\mathcal{R}_{{\Gamma}_{11}}(\mathcal{P}_2(u_2))|)=(\frac{9+3}{2},\frac{9-3}{2})=(6,3),$$
$$(|\mathcal{L}_{{\Gamma}_{11}}(\mathcal{P}_3(u_2))|,|\mathcal{R}_{{\Gamma}_{11}}(\mathcal{P}_3(u_2))|)=(\frac{9+1}{2},\frac{9-1}{2})=(5,4),$$
$$(|\mathcal{L}_{{\Gamma}_{11}}(\mathcal{P}_4(u_2))|,|\mathcal{R}_{{\Gamma}_{11}}(\mathcal{P}_4(u_2))|)=(\frac{9+1}{2},\frac{9+1}{2})=(5,5).$$
Then we balance the number of left arcs and the number of right arcs with respect to each vertex of $\mathcal{P}_1(u_1)$, $\mathcal{P}_2(u_1)$, $\mathcal{P}_3(u_1)$, $\mathcal{P}_4(u_1)$,
$\mathcal{P}_1(u_2)$, $\mathcal{P}_2(u_2)$,
$\mathcal{P}_3(u_2)$, $\mathcal{P}_4(u_2)$
by distorting seven edges (see Figure \ref{fig:gamma11togamma*v}), $$\mathcal{P}_1(u_1)\mathcal{P}_2(u_1), \mathcal{P}_2(u_1)\mathcal{P}_3(u_1), \mathcal{P}_3(u_1)\mathcal{P}_4(u_1),$$
$$\mathcal{P}_1(u_2)\mathcal{P}_2(u_2), \mathcal{P}_2(u_2)\mathcal{P}_3(u_2), \mathcal{P}_3(u_2)\mathcal{P}_4(u_2),$$
$$\mathcal{P}_4(u_1)\mathcal{P}_1(u_2),$$
`locally', and check that in ${\Gamma}_{11}^*$, $$(|\mathcal{L}_{{\Gamma}_{11}^*}(\mathcal{P}_1(u_1))|,|\mathcal{R}_{{\Gamma}_{11}^*}(\mathcal{P}_1(u_1))|)=(\frac{9+1}{2},\frac{9+3}{2})=(5,6),$$
$$(|\mathcal{L}_{{\Gamma}_{11}^*}(\mathcal{P}_2(u_1))|,|\mathcal{R}_{{\Gamma}_{11}^*}(\mathcal{P}_2(u_1))|)=(\frac{9+1}{2},\frac{9+3}{2})=(5,6),$$
$$(|\mathcal{L}_{{\Gamma}_{11}^*}(\mathcal{P}_3(u_1))|,|\mathcal{R}_{{\Gamma}_{11}^*}(\mathcal{P}_3(u_1))|)=(\frac{9+3}{2},\frac{9+1}{2})=(6,5),$$
$$(|\mathcal{L}_{{\Gamma}_{11}^*}(\mathcal{P}_4(u_1))|,|\mathcal{R}_{{\Gamma}_{11}^*}(\mathcal{P}_4(u_1))|)=(\frac{9+1}{2},\frac{9+3}{2})=(5,6),$$
and
$$(|\mathcal{L}_{{\Gamma}_{11}^*}(\mathcal{P}_1(u_2))|,|\mathcal{R}_{{\Gamma}_{11}^*}(\mathcal{P}_1(u_2))|)=(\frac{9+1}{2},\frac{9+3}{2})=(5,6),$$
$$(|\mathcal{L}_{{\Gamma}_{11}^*}(\mathcal{P}_2(u_2))|,|\mathcal{R}_{{\Gamma}_{11}^*}(\mathcal{P}_2(u_2))|)=(\frac{9+3}{2},\frac{9+1}{2})=(6,5),$$
$$(|\mathcal{L}_{{\Gamma}_{11}^*}(\mathcal{P}_3(u_2))|,|\mathcal{R}_{{\Gamma}_{11}^*}(\mathcal{P}_3(u_2))|)=(\frac{9+1}{2},\frac{9+3}{2})=(5,6),$$
$$(|\mathcal{L}_{{\Gamma}_{11}^*}(\mathcal{P}_4(u_2))|,|\mathcal{R}_{{\Gamma}_{11}^*}(\mathcal{P}_4(u_2))|)=(\frac{9+1}{2},\frac{9+3}{2})=(5,6).$$
\begin{figure}[!ht]
\centering\includegraphics[scale=1.0]{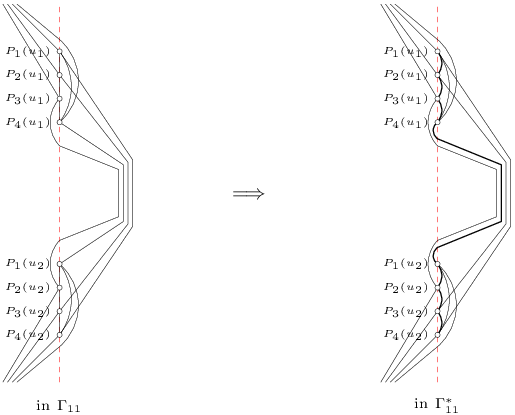} \caption{\small{Schematic diagram for the adjustments from $\Gamma_{11}$ to $\Gamma_{11}^*$ }}\label{fig:gamma11togamma*v}
\end{figure}

We see that the adjusted edges above remain self-symmetric.
Since $\mathcal{M}$ is a perfect matching and $u_1u_2$ is an arbitrary edge of the matching $\mathcal{M}$, we have that the number of left arcs and the number of right arcs with respect to each vertex in $\Gamma_{11}^*$ have been balanced,
and moreover, the above adjustments do not affect the number of fundamental structures in $\Gamma_{11}^*$ which is still the same as in $\Gamma_{11}$, i.e., the obtained drawing $\Gamma_{11}^*$ from $\Gamma_{11}$ satisfies the required Properties 1, 2.

Let $n\geq 13$. It remains to introduce the general rule on how to construct $\Gamma_n^*$ from $\Gamma_n$.
Let $u$ be an arbitrary vertex in the drawing $\Gamma_{n-2}^*$. Note that the mesh-like structure around $u^{(00)}u^{(10)}u^{(11)}u^{(01)}$ in $\widetilde{\Gamma}_n$ is depicted as Figure \ref{fig:mesh-likeStructurefromn-2ton}.
\begin{figure}[!ht]
\centering\includegraphics[scale=1.0]{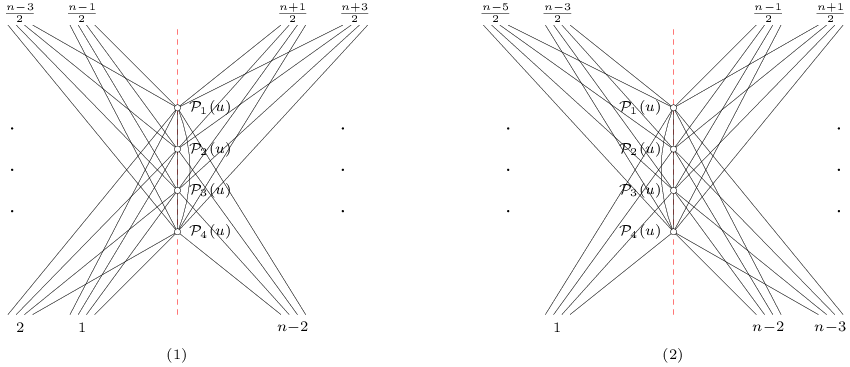} \caption{\small{The drawing of Mesh-like structure around $u^{(00)}u^{(10)}u^{(11)}u^{(01)}$ in $\widetilde{\Gamma}_n$ constructed from $\Gamma_{n-2}^*$}}\label{fig:mesh-likeStructurefromn-2ton}
\end{figure}
Since only the second kind of adjustments may be taken on edges in $\widetilde{\Gamma}_{n}$ (without the first kind of adjustments), due to the same reason as above for the rule of obtaining $\Gamma_{11}^*$ from  $\Gamma_{11}$, the number of left arcs and the number of right arcs with respect to each vertex of $u^{(00)},u^{(10)},u^{(11)},u^{(01)}$ in $\Gamma_n$ is the same as ones in $\widetilde{\Gamma}_n$, and moreover, the  small $4$-cycles $u^{(00)}u^{(10)}u^{(11)}u^{(01)}$ are not affected by the second kind of adjustment.  For the cases when the mesh-like structure around $u^{(00)}u^{(10)}u^{(11)}u^{(01)}$ in $\widetilde{\Gamma}_{n}$ is given as Diagram (1) or Diagram (2) in Figure \ref{fig:mesh-likeStructurefromn-2ton}, we check that in $\Gamma_{n}$, $$(|\mathcal{L}_{{\Gamma}_{n}}(\mathcal{P}_1(u))|,|\mathcal{R}_{{\Gamma}_{n}}(\mathcal{P}_1(u))|)=(\frac{n-1}{2},\frac{n-1}{2}),$$
$$(|\mathcal{L}_{{\Gamma}_{n}}(\mathcal{P}_2(u))|,|\mathcal{R}_{{\Gamma}_{n}}(\mathcal{P}_2(u))|)=(\frac{n-1}{2},\frac{n-3}{2}),$$
$$(|\mathcal{L}_{{\Gamma}_{n}}(\mathcal{P}_3(u))|,|\mathcal{R}_{{\Gamma}_{n}}(\mathcal{P}_3(u))|)=(\frac{n-1}{2},\frac{n-3}{2}),$$
$$(|\mathcal{L}_{{\Gamma}_{n}}(\mathcal{P}_4(u))|,|\mathcal{R}_{{\Gamma}_{n}}(\mathcal{P}_4(u))|)=(\frac{n-1}{2},\frac{n-1}{2}),$$
or
$$(|\mathcal{L}_{{\Gamma}_{n}}(\mathcal{P}_1(u))|,|\mathcal{R}_{{\Gamma}_{n}}(\mathcal{P}_1(u))|)=(\frac{n-1}{2},\frac{n-1}{2}),$$
$$(|\mathcal{L}_{{\Gamma}_{n}}(\mathcal{P}_2(u))|,|\mathcal{R}_{{\Gamma}_{n}}(\mathcal{P}_2(u))|)=(\frac{n-3}{2},\frac{n-1}{2}),$$
$$(|\mathcal{L}_{{\Gamma}_{n}}(\mathcal{P}_3(u))|,|\mathcal{R}_{{\Gamma}_{n}}(\mathcal{P}_3(u))|)=(\frac{n-3}{2},\frac{n-1}{2}),$$
$$(|\mathcal{L}_{{\Gamma}_{n}}(\mathcal{P}_4(u))|,|\mathcal{R}_{{\Gamma}_{n}}(\mathcal{P}_4(u))|)=(\frac{n-1}{2},\frac{n-1}{2}),$$
respectively. Then we balance the number of left arcs and the number of right arcs with respect to  each vertex of $u^{(00)},u^{(10)},u^{(11)},u^{(01)}$
by distorting three edges (see Figure \ref{fig:gamma-to-gamma*ngeq13}), $$\mathcal{P}_1(u)\mathcal{P}_2(u), \mathcal{P}_2(u)\mathcal{P}_3(u), \mathcal{P}_3(u)\mathcal{P}_4(u),$$
\begin{figure}[!ht]
\centering\includegraphics[scale=1.0]{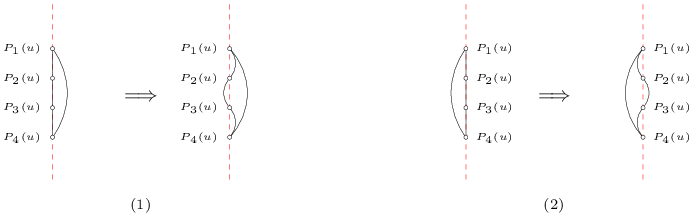} \caption{\small{Schematic diagrams for the adjustments from $\Gamma_n$ to $\Gamma_n^*$ with $n\geq 13$}}\label{fig:gamma-to-gamma*ngeq13}
\end{figure}
\\
and check that in ${\Gamma}_{n}^*$, $$(|\mathcal{L}_{{\Gamma}_{n}^*}(\mathcal{P}_1(u))|,|\mathcal{R}_{{\Gamma}_{n}^*}(\mathcal{P}_1(u))|)=(\frac{n-1}{2},\frac{n+1}{2}),$$
$$(|\mathcal{L}_{{\Gamma}_{n}^*}(\mathcal{P}_2(u))|,|\mathcal{R}_{{\Gamma}_{n}^*}(\mathcal{P}_2(u))|)=(\frac{n+1}{2},\frac{n-1}{2}),$$
$$(|\mathcal{L}_{{\Gamma}_{n}^*}(\mathcal{P}_3(u))|,|\mathcal{R}_{{\Gamma}_{n}^*}(\mathcal{P}_3(u))|)=(\frac{n+1}{2},\frac{n-1}{2}),$$
$$(|\mathcal{L}_{{\Gamma}_{n}^*}(\mathcal{P}_4(u))|,|\mathcal{R}_{{\Gamma}_{n}^*}(\mathcal{P}_4(u))|)=(\frac{n-1}{2},\frac{n+1}{2}),$$
or
$$(|\mathcal{L}_{{\Gamma}_{n}^*}(\mathcal{P}_1(u))|,|\mathcal{R}_{{\Gamma}_{n}^*}(\mathcal{P}_1(u))|)=(\frac{n+1}{2},\frac{n-1}{2}),$$
$$(|\mathcal{L}_{{\Gamma}_{n}^*}(\mathcal{P}_2(u))|,|\mathcal{R}_{{\Gamma}_{n}^*}(\mathcal{P}_2(u))|)=(\frac{n-1}{2},\frac{n+1}{2}),$$
$$(|\mathcal{L}_{{\Gamma}_{n}^*}(\mathcal{P}_3(u))|,|\mathcal{R}_{{\Gamma}_{n}^*}(\mathcal{P}_3(u))|)=(\frac{n-1}{2},\frac{n+1}{2}),$$
$$(|\mathcal{L}_{{\Gamma}_{n}^*}(\mathcal{P}_4(u))|,|\mathcal{R}_{{\Gamma}_{n}^*}(\mathcal{P}_4(u))|)=(\frac{n+1}{2},\frac{n-1}{2}),$$
for each case respectively. Similarly as above, we have that each edge in ${\Gamma}_{n}^*$ keeps self-symmetric and the drawing ${\Gamma}_{n}^*$ satisfies Properties 1, 2 as  desired.

This completes the whole inductive process of constructing $\Gamma_{n+2}$ from
$\Gamma_n$ for all odd $n\geq 5$.

\subsection{Construction of the drawing $\Gamma_{n+1}$ out of $\Gamma_n^*$ for all odd $n\geq 5$}

Let $n\geq 5$ be an odd integer. We shall construct the desired drawing $\Gamma_{n+1}$ directly
from the drawing $\Gamma_n^*$ given in Subsection 3.1. To make the process clear, we give a part of the drawing $\Gamma_{6}$ in Figure \ref{fig:aid diagramsGamma6} and the drawing $\Gamma_{8}$ in Figure \ref{fig:Gamma8} to illustrate the process of constructing $\Gamma_{6}$ from $\Gamma_5^*$ (see the corresponding Figure \ref{fig:aid diagrams1} (2)), and the process of constructing $\Gamma_{8}$ from $\Gamma_7^*$ (see the corresponding Figure \ref{fig:aid diagrams4}), respectively.

In general, the process of constructing $\Gamma_{n+1}$ from the drawing $\Gamma_n^*$ is as follows.

Let $u$ be an arbitrary vertex in the drawing $\Gamma_n^*$. We
locate the two new vertices $u^{(0)}$ and $u^{(1)}$ in $\Gamma_{n+1}$
with
$$X_{u^{(0)}}=X_{u^{(1)}}=X_u,$$
$$Y_{u^{(0)}}=Y_u,$$
$$Y_{u^{(1)}}=Y_{u}+\frac{Y_{\widehat{u}}-Y_{u}}{\mathcal{N}},$$
and the new edge $u^{(0)}u^{(1)}$ drawn precisely at the line $x=X_u$.

Let $u_1u_2$ be an arbitrary edge in the drawing $\Gamma_{n}^*$.
We draw the two edges $u_1^{(0)}u_2^{(0)}$ and $u_1^{(1)}u_2^{(1)}$
in $\Gamma_{n+1}$ to be a `bunch' such that the bunch is along the
original route of $u_1u_2$ in $\Gamma_n^*$. In particular, the natures of each one, say $u_1^{(a)}u_2^{(a)}$ where $a\in \{0,1\}$, of the two edges $u_1^{(0)}u_2^{(0)}$, $u_1^{(1)}u_2^{(1)}$,
with respect to both its ends $u_1^{(a)}, u_2^{(a)}$ in $\Gamma_{n+1}$ are the same as the natures of $u_1u_2$ with respect to both its ends $u_1, u_2$ in $\Gamma_n^*$. Similarly as in Step 2 of Subsection 3.1, since the edge $u_1u_2$ is self-symmetric in $\Gamma_n^*$,
we have that the two edges $u_1^{(0)}u_2^{(0)}$, $u_1^{(1)}u_2^{(1)}$ do not cross each other, i.e., \begin{equation}\label{equation the bunch no crossings even}
\nu_{\Gamma_{n+1}}\big{(}\{u_1^{(0)}u_2^{(0)},\ u_1^{(1)}u_2^{(1)}\}
\big{)}=0.
\end{equation}
\begin{figure}[!ht]
\centering\includegraphics[scale=1.0]{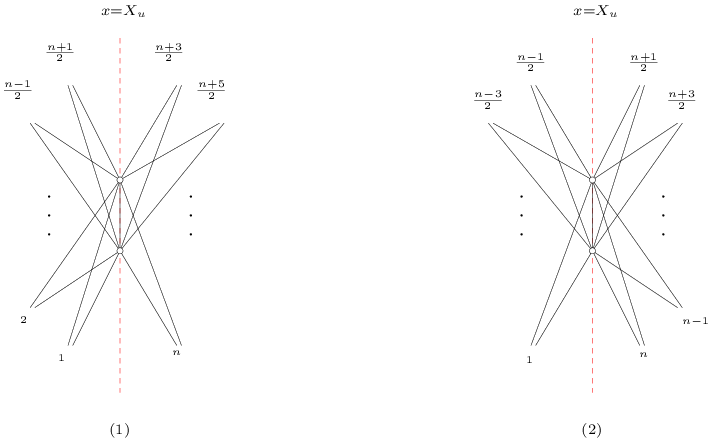} \caption{\small{The drawing of Mesh-like structure around $u^{(0)}u^{(1)}$ in $\Gamma_{n+1}$ constructed from $\Gamma_n^*$
}}\label{fig:Meshforeven}
\end{figure}
Moreover, we see that the mesh-like structure formed around $u^{(0)}u^{(1)}$ is depicted as Diagram (1) or Diagram (2) in Figure \ref{fig:Meshforeven} according to $(\mathcal{L}_{\Gamma_{n}^*}(u),\mathcal{R}_{\Gamma_{n}^*}(u))=(\frac{n+1}{2},\frac{n-1}{2})$ or $(\mathcal{L}_{\Gamma_{n}^*}(u),\mathcal{R}_{\Gamma_{n}^*}(u))=(\frac{n-1}{2},\frac{n+1}{2})$ respectively. This completes the description of the process in this subsection.

\begin{figure}
\centering
\includegraphics[scale=1.3]{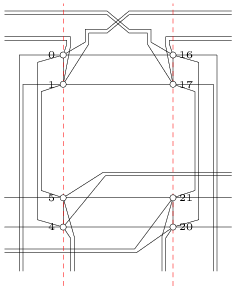}
\caption{\small{A part of $\Gamma_6$ illustrating the process to construct $\Gamma_6$ from $\Gamma_5^*$}}\label{fig:aid diagramsGamma6}
\end{figure}

\bigskip

\subsection{Calculations of the number of crossings in the desired drawing}

In this subsection, we shall calculate the number of crossings in
the drawing $\Gamma_n$ constructed in Subsection 3.1 and Subsection
3.2.

The following Lemma \ref{lemma crossings in MLk} can be found in \cite{FFSV08}. For the reader's convenience, we give a proof below.

\begin{lemma}\label{lemma crossings in MLk}
Let $n\geq 5$ be an odd integer.   Let $m_n$ be the number of crossings given in each diagram of Figure \ref{fig:Meshforeven}. Let $M_n$ be the number of crossings given in each diagram of Figure \ref{fig:local drawing}. Let $\widehat{M}_n$ be the number of crossings given in each diagram of Figure \ref{fig:basic-mesh-bold-adjusted}.
Then
$${\rm (i) } \ m_n={\frac{n+1}{2}\choose 2}+ {\frac{n-1}{2}\choose 2},$$
$${\rm (ii) } \ M_n={4 \choose 2}\cdot {\frac{n+1}{2}\choose 2}+{4\choose2}\cdot {\frac{n-1}{2}\choose 2}+(n-1)$$ and $${\rm (iii) } \  \widehat{M}_n=M_n-1.$$
\end{lemma}

\begin{proof} (i) \ We take Diagram (1) of Figure \ref{fig:Meshforeven} for example to show the calculations, because Diagram (2) is just a reflection of Diagram (1) and definitely has the same number of crossings.
Notice that any two distinct bunches which lies on the left side of the line $x=X_u$ (i.e., numbered from $1$ to $\frac{n+1}{2}$) have exactly one crossing. This implies that the number of crossings formed by the bunches on the left side of the line $x=X_u$ is $${\frac{n+1}{2}\choose 2}.$$ Similarly, the number of crossings formed by the bunches on the right side of the line $x=X_u$ (i.e., numbered from $\frac{n+3}{2}$ to $n$) is $${\frac{n-1}{2}\choose 2}.$$ Hence, this gives $m_n={\frac{n+1}{2}\choose 2}+ {\frac{n-1}{2}\choose 2}$.

(ii) \ We take Diagram (2) of Figure \ref{fig:local drawing} for example to show the calculations.
Notice that any two distinct bunches which are numbered from $1$ to $\frac{n-1}{2}$ have
exactly $6={4 \choose 2}$ crossings, and thus,  the number of crossings formed by the bunches (numbered from $1$ to $\frac{n-1}{2}$) is $${4\choose 2}\cdot {\frac{n-1}{2}\choose 2}.$$ Moreover, the arc $\mathcal{P}_1(u)\mathcal{P}_4(u)$ has exactly two crossings with each bunch numbered from $1$ to $\frac{n-1}{2}$, that is, the corresponding number of crossings is $$2* \frac{n-1}{2}=n-1.$$  In a similar observation, we have that the number of crossings formed by the bunches numbered from $\frac{n+1}{2}$ to $n$ is $${4\choose 2}\cdot {\frac{n+1}{2}\choose 2}.$$ Therefore, this gives $M_n={4 \choose 2}\cdot {\frac{n+1}{2}\choose 2}+{4\choose 2}\cdot {\frac{n-1}{2}\choose 2}+(n-1).$

\begin{figure}
\centering
\includegraphics[scale=1.0]{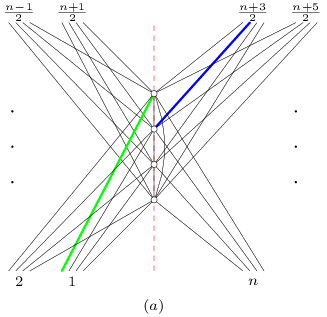} \hspace{50pt}
\includegraphics[scale=1.0]{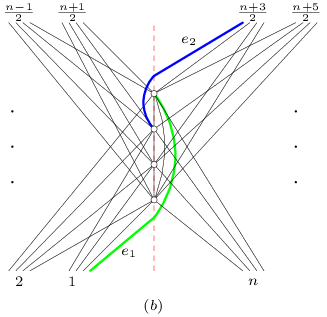}
\caption{\small{Auxiliary drawings for the calculation to show  $\widehat{M}_n=M_n-1$}}\label{fig:calculationforMnandMntilde}
\end{figure}

(iii) \ We take Diagram (3) of Figure \ref{fig:basic-mesh-bold-adjusted} for example show the calculations.  Notice that we obtain this diagram from Diagram (3) of Figure \ref{fig:basic-mesh-bold} by distorting two edges. For the notational convenience, we put both diagrams into a new figure (see Figure \ref{fig:calculationforMnandMntilde}) and marked the adjusted edges in colors.
Note that we distort only two edges, say $e_1$ and $e_2$ (emphasized in green and blue respectively) during the process from Diagram (a) to Diagram (b) in Figure \ref{fig:calculationforMnandMntilde}. To show Conclusion (iii), we need to calculate the crossings formed by $e_1$ and by $e_2$ in both  Diagram (a) and Diagram (b).

Observe Diagram (a).  The edge $e_1$ has exactly $3$ crossings with each bunch numbered from $2$ to $\frac{n+1}{2}$.
The edge $e_2$ has exactly one crossing with each bunch numbered from $\frac{n+5}{2}$ to $n$, and moreover, has one crossing with the arc $\mathcal{P}_1(u)\mathcal{P}_4(u)$. Hence, the total number of crossings formed by $e_1$ and by $e_2$ within the mesh is
\begin{equation}\label{equation crossings in the basic mesh}
3\cdot \frac{n-1}{2}+(\frac{n-3}{2}+1)=2n-2.
\end{equation}

Now observe Diagram (b).   The edge $e_1$ has exactly $3$ crossings with each bunch numbered from $\frac{n+5}{2}$ to $n$, and has $2$ crossings with the bunch numbered $\frac{n+3}{2}$ (crossed the two edges which lie in the bunch numbered $\frac{n+3}{2}$ and are incident to $\mathcal{P}_3(u)$ and $\mathcal{P}_4(u)$ respectively).
The edge $e_2$ has exactly one crossing with each bunch numbered from $2$ to $\frac{n+1}{2}$. Hence, the total number of crossings formed by $e_1$ and $e_2$ within the mesh is
\begin{equation}\label{equation crossings in the adjusted mesh}
(3\cdot \frac{n-3}{2}+2)+\frac{n-1}{2}=2n-3.
\end{equation}

Then Conclusion (iii) follows from \eqref{equation crossings in the basic mesh} and \eqref{equation crossings in the adjusted mesh} readily.
\end{proof}

To make the calculations clear,  we also give the drawings (see Figures \ref{fig:mesh7}, \ref{fig:mesh8}, \ref{fig:mesh9}, \ref{fig:mesh10}, \ref{fig:mesh11}, \ref{fig:mesh12}) and the corresponding number of crossings (see Table \ref{table:values of mesh}) of the mesh-like structures as depicted in Figures \ref{fig:local drawing} and \ref{fig:Meshforeven} with $n\in \{7,9,11\}$.

\begin{figure}
\centering
\includegraphics[scale=1.1]{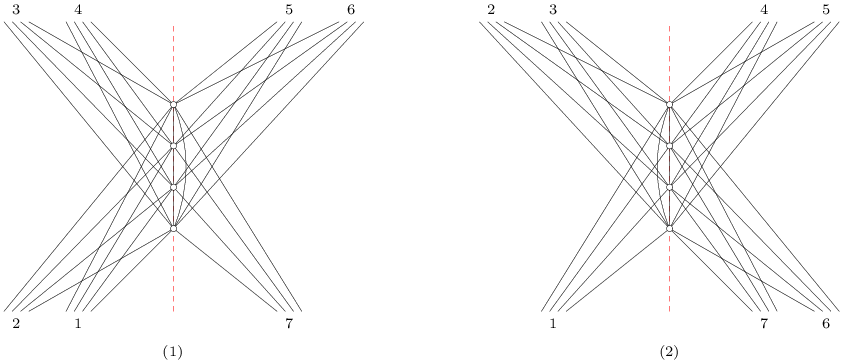}
\caption{\small{The drawing of Mesh-like structure around $u^{(0)}u^{(1)}u^{(2)}u^{(3)}$ in $\widetilde{\Gamma}_9$ constructed from $\Gamma_7^*$}}\label{fig:mesh7}
\end{figure}

\begin{figure}
\centering
\includegraphics[scale=1.0]{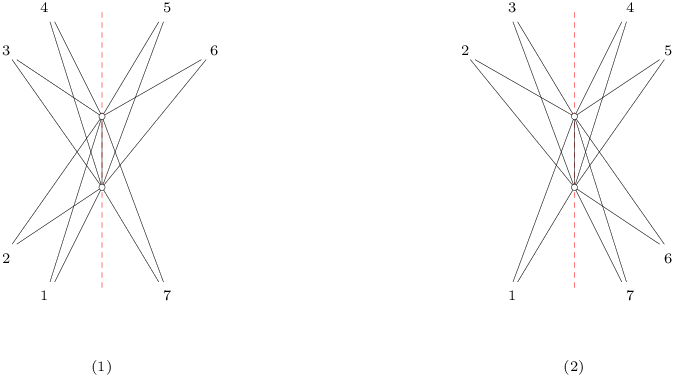}
\caption{\small{The drawing of Mesh-like structure around $u^{(0)}u^{(1)}$ in $\Gamma_8$ constructed from $\Gamma_7^*$}}\label{fig:mesh8}
\end{figure}

\begin{figure}
\centering
\includegraphics[scale=1.0]{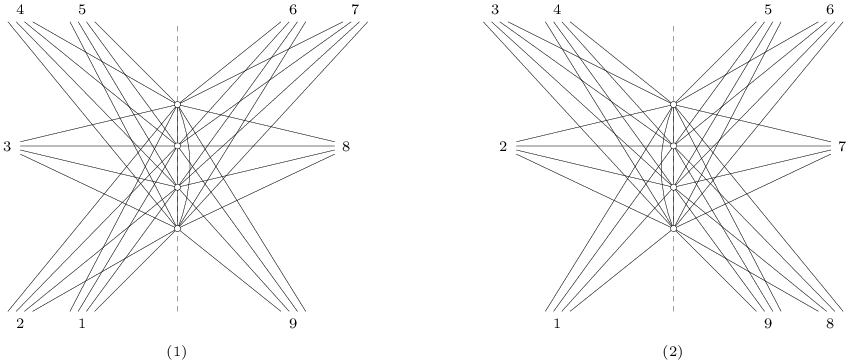}
\caption{\small{The drawing of Mesh-like structure around $u^{(0)}u^{(1)}u^{(2)}u^{(3)}$ in $\widetilde{\Gamma}_{11}$ constructed from $\Gamma_9^*$}}\label{fig:mesh9}
\end{figure}

\begin{figure}
\centering
\includegraphics[scale=1.0]{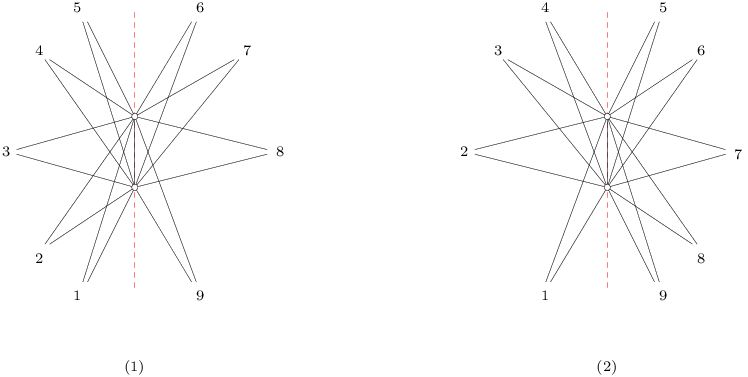}
\caption{\small{The drawing of Mesh-like structure around $u^{(0)}u^{(1)}$ in $\Gamma_{10}$ constructed from $\Gamma_9^*$}}\label{fig:mesh10}
\end{figure}

\begin{figure}
\centering
\includegraphics[scale=1.0]{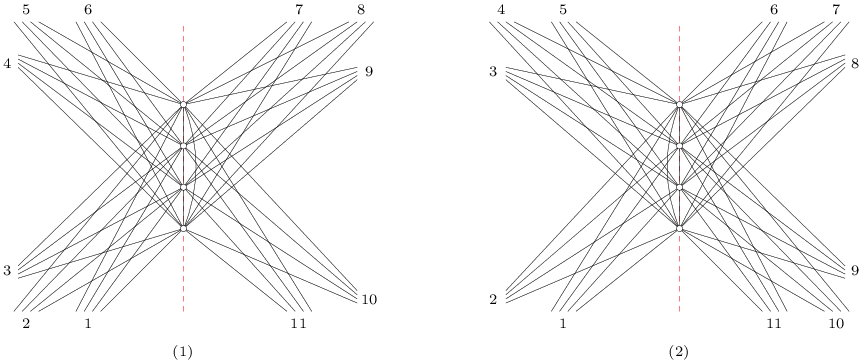}
\caption{\small{The drawing of Mesh-like structure around $u^{(0)}u^{(1)}u^{(2)}u^{(3)}$ in $\widetilde{\Gamma}_{13}$ constructed from $\Gamma_{11}^*$}}\label{fig:mesh11}
\end{figure}

\begin{figure}
\centering
\includegraphics[scale=1.0]{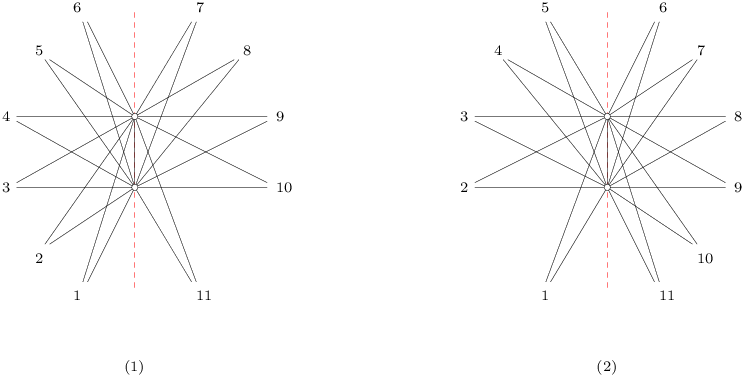}
\caption{\small{The drawing of Mesh-like structure around $u^{(0)}u^{(1)}$ in $\Gamma_{12}$ constructed from $\Gamma_{11}^*$}}\label{fig:mesh12}
\end{figure}

\medskip

\begin{table}[htbp]
\centering \small{
\begin{tabular}{|c|c|c|c|}
  \hline
  % after \\: \hline or \cline{col1-col2} \cline{col3-col4} ...
  $n$     & \ \ $m_n$  \ \ & \ \ \ \ \ \ $M_n$ \ \ \ \ \ \ \\ \hline
  $7$     & 9 & 60 \\ \hline
  $9$     & 16 & 104  \\ \hline
  $11$    & 25 & 160 \\ \hline
  \end{tabular}}
\caption{\small{Values of $m_n$ and $M_n$ for $n\in \{7,9,11\}$}}\label{table:values of mesh}
\end{table}

\bigskip

\bigskip

\bigskip

\begin{lemma}\label{Lemma calculation for odd n} For odd integers $n\geq 5$,  we define a sequence of numbers $A_n$ satisfying the recurrence relation and the initial condition as follows:

(i) \begin{equation}\label{equation An}
A_{n+2}=16\cdot A_n+ 2^n\cdot \big{(}{4
\choose 2}\cdot {\frac{n+1}{2}\choose 2}+{4 \choose 2}\cdot
{\frac{n-1}{2}\choose 2}+(n-1)\big{)}- \epsilon_n \cdot 2^{n}-8\cdot
2^{\frac{n-3}{2}}
\end{equation}
 where
\begin{equation}\label{equation definition of epsilon n}
\begin{array}{llll} \epsilon_n= &\left
\{\begin{array}{llll}
               1, & \mbox{ if } n\in \{5,7,9\}; \\
               0, & \mbox{ if } n\geq 11. \\
                       \end{array}
           \right. \\
\end{array}
\end{equation}

(ii) $A_5=56.$

Then we have $$\begin{array}{llll} & A_n=\left \{\begin{array}{llll}
               \frac{139}{896} 4^{n}-(\frac{n^2+1}{2}) 2^{n-2}+\frac{1}{7}\cdot 2^{\frac{n+1}{2}}, & \mbox{ if }\  n\in \{5,7,9,11\};\\
               \frac{139}{896} 4^{n}-(\frac{n^2+1}{2}) 2^{n-2}+\frac{1}{7}\cdot 2^{\frac{n+1}{2}}+\frac{2^{n-11}-1}{3}\cdot 2^{n-2}, &  \mbox{ otherwise.}\  \\
              \end{array}
           \right. \\
\end{array}$$.

\end{lemma}

\begin{proof} \ For odd integers $n\geq 5$,  we define a sequence of numbers $B_n$ satisfying the recurrence relation and the initial condition given by
\begin{equation}\label{equation definition of Bn}
B_{n+2}=16\cdot B_n+ 2^n\cdot \big{(}{4
\choose 2}\cdot {\frac{n+1}{2}\choose 2}+{4 \choose 2}\cdot
{\frac{n-1}{2}\choose 2}+(n-1)\big{)}-  2^{n}-8\cdot
2^{\frac{n-3}{2}}
\end{equation}
and
 \begin{equation}\label{equation B5=A5=56}
 B_5=A_5=56.
 \end{equation}

By \eqref{equation An}, \eqref{equation definition of epsilon n}, \eqref{equation definition of Bn} and \eqref{equation B5=A5=56}, we see that \begin{equation}\label{equation An=Bnforn=5to11}
A_n=B_n \mbox{ for } n\in \{5,7,9,11\}.
\end{equation}

 Now we show that
 \begin{equation}\label{equation Bn=}
 B_n=\frac{139}{896} 4^{n}-(\frac{n^2+1}{2}) 2^{n-2}+\frac{1}{7}\cdot 2^{\frac{n+1}{2}} \ \ \ \mbox{ for any odd integer } n\geq 5.
 \end{equation}
We shall prove \eqref{equation Bn=} by induction on $n$. If $n=5$, we verify that $B_5=\frac{139}{896} 4^{5}-(\frac{5^2+1}{2}) 2^{5-2}+\frac{1}{7}\cdot 2^{\frac{5+1}{2}}=56$, done.
Suppose $n\geq 5$ is an odd integer and \eqref{equation Bn=} holds for $n$. It suffices to prove \eqref{equation Bn=} holds for $n+2$.
By \eqref{equation definition of Bn}  and the hypothesis, we have that
$$\begin{array}{llll}
& &B_{n+2}
\\
&=& 16\cdot B_n+ 2^n\cdot \left({4
\choose 2}\cdot {\frac{n+1}{2}\choose 2}+{4 \choose 2}\cdot
{\frac{n-1}{2}\choose 2}+(n-1)\right)- 2^{n}-8\cdot
2^{\frac{n-3}{2}} \\
&=&
16\cdot \left( \frac{139}{896} 4^{n}-(\frac{n^2+1}{2}) \cdot 2^{n-2}+\frac{1}{7}\cdot 2^{\frac{n+1}{2}}\right)+ 2^n\cdot \left({4
\choose 2}\cdot {\frac{n+1}{2}\choose 2}+{4 \choose 2}\cdot
{\frac{n-1}{2}\choose 2}+(n-1)\right)- 2^{n}-8\cdot
2^{\frac{n-3}{2}} \\
&=&
\left(\frac{139}{896} 4^{n+2}-(2 n^2+2) \cdot 2^n+\frac{8}{7}\cdot 2^{\frac{n+3}{2}}\right )+ 2^n\cdot \big{(}\frac{3n^2-4n+1}{2}
\big{)}- 2^{n}-
2^{\frac{n+3}{2}} \\
&=&
\frac{139}{896} 4^{n+2}+\left(-(2 n^2+2)+\frac{3n^2-4n+1}{2}-1\right) \cdot 2^n+(\frac{8}{7}-1)\cdot 2^{\frac{n+3}{2}}\\
&=&
\frac{139}{896} 4^{n+2}-\left(\frac{(n+2)^2+1}{2}\right) \cdot 2^{n}+\frac{1}{7}\cdot 2^{\frac{n+3}{2}}\\
&=&
\frac{139}{896} 4^{n+2}-\left(\frac{(n+2)^2+1}{2}\right) \cdot 2^{(n+2)-2}+\frac{1}{7}\cdot 2^{\frac{(n+2)+1}{2}}, \\
\end{array}$$
which proves \eqref{equation Bn=}.

Next we show that
\begin{equation}\label{equation An=Bn+}
A_n=B_n+\frac{2^{n-11}-1}{3}\cdot 2^{n-2} \ \ \ \mbox{ for any odd integer } n\geq 13.
\end{equation}
By induction on $n$. If $n=13$, it follows from  \eqref{equation An}, \eqref{equation definition of epsilon n},
\eqref{equation definition of Bn} and \eqref{equation An=Bnforn=5to11} that
$$\begin{array}{llll}
A_n&=&A_{13}
\\
&=& 16\cdot A_{11}+ 2^{11}\cdot \big{(}{4
\choose 2}\cdot {\frac{11+1}{2}\choose 2}+{4 \choose 2}\cdot
{\frac{11-1}{2}\choose 2}+(11-1)\big{)}-8\cdot
2^{\frac{11-3}{2}} \\
&=&
16\cdot B_{11}+ 2^{11}\cdot \big{(}{4
\choose 2}\cdot {\frac{11+1}{2}\choose 2}+{4 \choose 2}\cdot
{\frac{11-1}{2}\choose 2}+(11-1)\big{)}-8\cdot
2^{\frac{11-3}{2}} \\
&=&
\left(16\cdot B_{11}+ 2^{11}\cdot \big{(}{4
\choose 2}\cdot {\frac{11+1}{2}\choose 2}+{4 \choose 2}\cdot
{\frac{11-1}{2}\choose 2}+(11-1)\big{)}-2^{11}-8\cdot
2^{\frac{11-3}{2}}\right)+2^{11} \\
&=&
B_{13}+2^{11} \\
&=&
B_n+\frac{2^{n-11}-1}{3}\cdot 2^{n-2},
\end{array}$$ done.
Suppose $n\geq 13$ is an odd integer and \eqref{equation An=Bn+} holds for $n$. It suffices to prove \eqref{equation An=Bn+} holds for $n+2$. By  \eqref{equation An}, \eqref{equation definition of epsilon n},
\eqref{equation definition of Bn} and the hypothesis, we conclude that
$$\begin{array}{llll}
& &A_{n+2}
\\
&=& 16\cdot A_n+ 2^n\cdot \big{(}{4
\choose 2}\cdot {\frac{n+1}{2}\choose 2}+{4 \choose 2}\cdot
{\frac{n-1}{2}\choose 2}+(n-1)\big{)}-8\cdot
2^{\frac{n-3}{2}} \\
&=&
16\cdot \left( B_n+\frac{2^{n-11}-1}{3}\cdot 2^{n-2}\right)+ 2^n\cdot \left({4
\choose 2}\cdot {\frac{n+1}{2}\choose 2}+{4 \choose 2}\cdot
{\frac{n-1}{2}\choose 2}+(n-1)\right)-8\cdot
2^{\frac{n-3}{2}} \\
&=&
16\cdot B_n+16\cdot \frac{2^{n-11}-1}{3}\cdot 2^{n-2}+ 2^n\cdot \left({4
\choose 2}\cdot {\frac{n+1}{2}\choose 2}+{4 \choose 2}\cdot
{\frac{n-1}{2}\choose 2}+(n-1)\right)-8\cdot
2^{\frac{n-3}{2}} \\
&=&
16\cdot B_n+\frac{2^{(n+2)-11}-4}{3}\cdot 2^n+ 2^n\cdot \left({4
\choose 2}\cdot {\frac{n+1}{2}\choose 2}+{4 \choose 2}\cdot
{\frac{n-1}{2}\choose 2}+(n-1)\right)-8\cdot
2^{\frac{n-3}{2}} \\
&=&
16\cdot B_n+\left(\frac{2^{(n+2)-11}-1}{3}\cdot 2^{n}-2^n\right)+ 2^n\cdot \left({4
\choose 2}\cdot {\frac{n+1}{2}\choose 2}+{4 \choose 2}\cdot
{\frac{n-1}{2}\choose 2}+(n-1)\right)-8\cdot
2^{\frac{n-3}{2}} \\
&=&
16\cdot B_n+2^n\cdot \left({4
\choose 2}\cdot {\frac{n+1}{2}\choose 2}+{4 \choose 2}\cdot
{\frac{n-1}{2}\choose 2}+(n-1)\right)-2^n-8\cdot
2^{\frac{n-3}{2}}+\frac{2^{(n+2)-11}-1}{3}\cdot 2^{n} \\
&=& B_{n+2}+\frac{2^{(n+2)-11}-1}{3}\cdot 2^{n}\\
&=& B_{n+2}+\frac{2^{(n+2)-11}-1}{3}\cdot 2^{(n+2)-2},
\end{array}$$
which proves \eqref{equation An=Bn+}. \qed

Then  the lemma follows from \eqref{equation An=Bnforn=5to11}, \eqref{equation Bn=} and \eqref{equation An=Bn+} readily.
\end{proof}

\medskip

Now we are in a position to give the detailed calculation of $\nu(\Gamma_n)$.
By the process described in Subsection 3.1, we conclude that for all odd
integers $n\geq 5$,
\begin{equation}\label{equation gamma=gamma*}
\nu(\Gamma_n^*)=\nu(\Gamma_n)
\end{equation}
and
\begin{equation}\label{equation crossingsgamman+2}
\nu(\widetilde{\Gamma}_{n+2})=16\cdot \nu(\Gamma_n^*)+ 2^n \cdot M_n.
\end{equation}
where $16=4^2$ in the term $16\cdot \nu(\Gamma_n^*)$ is the crossings produced by any two bunches, say $$u_1^{(00)}u_2^{(00)},u_1^{(10)}u_2^{(10)},u_1^{(11)}u_2^{(11)},u_1^{(01)}u_2^{(01)}$$
and $$u_3^{(00)}u_4^{(00)},u_3^{(00)}u_4^{(00)},u_3^{(00)}u_4^{(00)},u_3^{(00)}u_4^{(00)}$$
with $u_1u_2$ and $u_3u_4$ cross in $\Gamma_n^*$,
and where the term $2^n \cdot M_n$ is the total crossings in the mesh-like structure around  $u^{(00)}u^{(10)}u^{(11)}u^{(01)}$ for all vertices $u$ in $\Gamma_n^*$.

By Property 2 and Property 3 holding for $\Gamma_n^*$, since we applied the first kind of adjustments for $\widetilde{\Gamma}_{n+2}$ when $n\in \{5,7,9\}$ and applied the second kind of adjustments for $\widetilde{\Gamma}_{n+2}$ with all odd $n\geq 5$, it follows from  Conclusion {\rm (iii)} of Lemma \ref{lemma crossings in MLk} and \eqref{equation crossings inner and outer}, \eqref{equation crossings inner crossings same} that $$\nu(\Gamma_{n+2})=\nu(\widetilde{\Gamma}_{n+2})-\epsilon_n\cdot
2^{n}-8\cdot 2^{\frac{n-3}{2}},$$ where
$$\begin{array}{llll} \epsilon_n= &\left
\{\begin{array}{llll}
               1, & \mbox{ if } n\in \{5,7,9\}; \\
               0, & \mbox{ if } n\geq 11, \\
                       \end{array}
           \right. \\
\end{array}$$
Combined with \eqref{equation gamma=gamma*}, \eqref{equation crossingsgamman+2} and Conclusion {\rm (ii)} of Lemma \ref{lemma crossings in MLk}, we have
\begin{eqnarray}\label{equation computations for n=odd}
\nu(\Gamma_{n+2})&=&16\cdot \nu(\Gamma_n)+ 2^n\cdot \big{(}{4
\choose 2}\cdot {\frac{n+1}{2}\choose 2}+{4 \choose 2}\cdot
{\frac{n-1}{2}\choose 2}+(n-1)\big{)}-\epsilon_n\cdot 2^{n}-8\cdot
2^{\frac{n-3}{2}}.\nonumber \\
\end{eqnarray}

Similarly as above, by \eqref{equation the bunch no crossings even}, Conclusion {\rm (i)} of Lemma \ref{lemma crossings in MLk} and the process described in Subsection 3.2, we can derive
that for all odd number $n\geq 5$,
\begin{eqnarray}\label{even crossings in both gamma and widetilde}
\nu(\Gamma_{n+1})&=&4\cdot \nu(\Gamma_{n}^*)+ 2^{n}\cdot
\big{(}{(n+1)/2\choose 2}+{(n-1)/2\choose 2}\big{)}\nonumber\\
&=&4\cdot \nu(\Gamma_{n})+2^{n-2}\cdot (n-1)^2.
\end{eqnarray}

Set $$\begin{array}{llll} \lambda_n= &\left
\{\begin{array}{llll}
               0, & \mbox{ if } 5\leq n\leq 12; \\
               1, & \mbox{ if } n\geq 13. \\
                       \end{array}
           \right. \\
\end{array}$$

By Lemma \ref{Lemma calculation for odd n}, we have that for  any odd integer $n\geq 5$,
 \begin{equation}\label{equation value for n is odd}
 \begin{array}{llll}
& &\nu(\Gamma_n)
\\
&=& \frac{139}{896} 4^{n}-(\frac{n^2+1}{2}) 2^{n-2}+\frac{1}{7}\cdot 2^{\frac{n+1}{2}}+\lambda_n\cdot \frac{2^{n-11}-1}{3}\cdot 2^{n-2} \\
&=& \frac{139}{896} 4^{n}-(\frac{n^2+1}{2}) 2^{n-2}+\frac{4}{7}\cdot 2^{\frac{n-3}{2}}+\lambda_n\cdot (\frac{2^{n-11}\cdot 2^{n-2}}{3}-\frac{2^{n-2}}{3}) \\
&=& \frac{139}{896} 4^{n}-(\frac{n^2+1}{2}) 2^{n-2}+\frac{4}{7}\cdot 2^{3\cdot (\frac{n-1}{2}) -n}+\lambda_n\cdot (\frac{2^n\cdot 2^{n}}{2^{13}\cdot 3}-\frac{2^{n-1}}{6}) \\
&=&
\frac{139}{896} 4^{n}-\lfloor\frac{n^2+1}{2}\rfloor 2^{n-2}+\frac{4}{7}\cdot 2^{3\lfloor\frac{n}{2}\rfloor -n}+\lambda_n\cdot (\frac{4^n}{24576} -\frac{4^{\lfloor\frac{ n}{2}\rfloor}}{6}). \\
\end{array}
\end{equation}

Notice that $$\lambda_{n-1}=\lambda_n \mbox{ for any even integer } n\geq 6.$$
Combined with \eqref{even crossings in both gamma and widetilde} and \eqref{equation value for n is odd}, we have that for any even integer $n\geq 6$,
\begin{equation}\label{equation value for n is even}
 \begin{array}{llll}
& &\nu(\Gamma_n)
\\
&=& 4\cdot \nu(\Gamma_{n-1})+2^{n-3}\cdot (n-2)^2 \\
&=& 4\cdot \left(\frac{139}{896} 4^{n-1}-\lfloor\frac{(n-1)^2+1}{2}\rfloor 2^{n-3}+\frac{4}{7}\cdot 2^{3\lfloor\frac{n-1}{2}\rfloor -(n-1)}+\lambda_{n-1}(\frac{4^{n-1}}{24576} -\frac{4^{\lfloor\frac{ n-1}{2}\rfloor}}{6})\right)+2^{n-3}\cdot (n-2)^2\\
&=& 4\cdot \left(\frac{139}{896} 4^{n-1}-\lfloor\frac{n^2-2n+2}{2}\rfloor 2^{n-3}+\frac{4}{7}\cdot 2^{3\cdot \frac{n-2}{2} -(n-1)}+\lambda_{n}(\frac{4^{n-1}}{24576} -\frac{4^{\frac{ n-2}{2}}}{6})\right)+(n^2-4n+4)\cdot 2^{n-3} \\
&=& 4\cdot \left(\frac{139}{896} 4^{n-1}-\frac{n^2-2n+2}{2}\cdot  2^{n-3}+\frac{4}{7}\cdot 2^{\frac{n-4}{2}}+\lambda_{n}(\frac{4^{n-1}}{24576} -\frac{4^{\frac{ n-2}{2}}}{6})\right)+\frac{n^2-4n+4}{2}\cdot 2^{n-2} \\
&=& \frac{139}{896} 4^n-(n^2-2n+2)\cdot  2^{n-2}+\frac{4}{7}\cdot 2^{\frac{n}{2}}+\lambda_{n}(\frac{4^{n}}{24576} -\frac{4^{\frac{ n}{2}}}{6})+\frac{n^2-4n+4}{2}\cdot 2^{n-2} \\
&=& \frac{139}{896} 4^n+\left(\frac{n^2-4n+4}{2}-(n^2-2n+2)\right)\cdot  2^{n-2}+\frac{4}{7}\cdot 2^{\frac{n}{2}}+\lambda_{n}(\frac{4^{n}}{24576} -\frac{4^{\frac{ n}{2}}}{6}) \\
&=& \frac{139}{896} 4^n-\frac{n^2}{2}\cdot 2^{n-2}+\frac{4}{7}\cdot 2^{3\cdot \frac{n}{2} -n}+\lambda_{n}(\frac{4^{n}}{24576} -\frac{4^{\frac{ n}{2}}}{6}) \\
&=& \frac{139}{896} 4^n-\lfloor\frac{n^2+1}{2}\rfloor\cdot 2^{n-2}+\frac{4}{7}\cdot 2^{3\lfloor\frac{n}{2}\rfloor -n}+\lambda_{n}(\frac{4^{n}}{24576} -\frac{4^{\lfloor\frac{ n}{2}}\rfloor}{6}). \\
\end{array}
\end{equation}
Combined  \eqref{equation value for n is odd} and \eqref{equation value for n is even},  we conclude that for all integers $n\geq 5$, $${\rm cr}(Q_n)\leq\nu(\Gamma_n)=\frac{139}{896} 4^{n}-\lfloor\frac{n^2+1}{2}\rfloor 2^{n-2}+\frac{4}{7}\cdot 2^{3\lfloor\frac{n}{2}\rfloor -n}+\lambda_n(\frac{4^n}{24576} -\frac{4^{\lfloor\frac{ n}{2}\rfloor}}{6}),$$
completing the calculations.

\section{Concluding remarks}

We first remark that the drawing for $Q_n$ constructed in this manuscript is not optimal for $n\geq 11$, that is, the upper
bound given in Theorem \ref{Theorem main theorem} can be improved
slightly by still applying that kind of adjustments in Step 3 for the progress of obtaining $\Gamma_{n+2}$ from
$\widetilde{\Gamma}_{n+2}$ for odd integers $n\geq 11$. In fact, for
$n\geq 11$  we still can find several disjoint enclosed
cycles, say $\mathcal{C}_1,\ldots,\mathcal{C}_m$ in $\Gamma_n^*$, unfortunately,
with $V(\mathcal{C}_1)\cup\cdots\cup V(\mathcal{C}_m)\subsetneq
V(Q_n)$. Theoretically, that would be insignificant since we have no
general rule for finding those enclosed cycles for all $n\geq 11$. Indeed,
the idea ``\emph{enclosed cycles}'' essentially is just an improved
variant of the method employed in \cite{FFSV08}.
The idea employed by L. Faria, C.M.H. de Figueiredo, O. S\'{y}kora, I. Vr\v{t}o in \cite{FFSV08} enlightened us to obtain the present drawing of $Q_n$ in this manuscript. So, the most important contribution on the Erd\H{o}s and Guy's problem was made by the four authors in \cite{FFSV08}. We remark that our drawing managed to use the idea of  `a fundamental structure' to decrease the crossings furthermore.

Note that the drawing in this manuscript is better than the conjectured crossings by Erd\H{o}s and Guy. Observed that the two coefficients $\frac{139}{896}$ (for the case $5\leq n\leq 12$) and $\frac{139}{896}+\frac{1}{24576}=\frac{26695}{172032}$ (for the case $n\geq 13$) of the leading term  $4^n$ in Theorem \ref{Theorem main theorem} are all less than that  coefficient $\frac{5}{32}=\frac{140}{896}=\frac{26880}{172032}$ conjectured by Erd\H{o}s and Guy.
In particular, denote $f(n)=\frac{5}{32}4^n-\lfloor\frac{n^2+1}{2}\rfloor 2^{n-2}$ to be the function of the conjectured values, and denote $$\Delta(n)=f(n)-\nu(\Gamma_n)$$ to be the difference function of the conjectured values and the number of crossings in our constructed drawing, that is,
$$\begin{array}{llll}
\Delta(n)&=& f(n)-\nu(\Gamma_n) \\
&=& \left(\frac{5}{32}4^n-\lfloor\frac{n^2+1}{2}\rfloor 2^{n-2}\right)-\left(\frac{139}{896} 4^{n}-\lfloor\frac{n^2+1}{2}\rfloor 2^{n-2}+\frac{4}{7}\cdot 2^{3\lfloor\frac{n}{2}\rfloor -n}+\lambda_n(\frac{4^n}{24576} -\frac{4^{\lfloor\frac{ n}{2}\rfloor}}{6})\right) \\
&=& \left(\frac{5}{32}-\frac{139}{896}-\lambda_n\frac{1}{24576}\right)4^n+(\lambda_n \frac{4^{\lfloor\frac{ n}{2}\rfloor}}{6}-\frac{4}{7}\cdot 2^{3\lfloor\frac{n}{2}\rfloor -n}) \\
\end{array}$$
For $n\geq 13$, since $\lambda_n=1$ , it follows that
$$\begin{array}{llll}
\Delta(n)&=& \left(\frac{5}{32}-\frac{139}{896}-\frac{1}{24576}\right)4^n+\left(\frac{4^{\lfloor\frac{ n}{2}\rfloor}}{6}-\frac{4}{7}\cdot 2^{3\lfloor\frac{n}{2}\rfloor -n}\right) \\
&=& \left(\frac{26880}{172032}-\frac{26688}{172032}-\frac{7}{172032}\right)4^n+\left(\frac{4^{\lfloor\frac{ n}{2}\rfloor}}{6}-\frac{4}{7}\cdot 2^{2\lfloor\frac{n}{2}\rfloor}\cdot 2^{\lfloor\frac{n}{2}\rfloor -n}\right) \\
&=& \frac{185}{172032}4^n+4^{\lfloor\frac{n}{2}\rfloor}\left(\frac{1}{6}-\frac{4}{7}\cdot 2^{\lfloor\frac{n}{2}\rfloor -n}\right) \\
&>& \frac{185}{172032}4^n. \\
\end{array}$$
Moreover, we conclude that $$\Delta(n)=\frac{185}{172032}4^n+{\rm O}(2^n) \ \ \ (n\rightarrow \infty).$$

Finally, we close this paper by a table to show the difference between
the number of crossings in our drawing and the values conjectured by
Erd\H{o}s and Guy for $5\leq n\leq 13$. \\

\begin{table}[htbp]
\centering \small{
\begin{tabular}{|c|c|c|c|}
  \hline
  % after \\: \hline or \cline{col1-col2} \cline{col3-col4} ...
  $n$     & \ \ Conjectured values  \ \ & \ \ \ \ \ \ Our results \ \ \ \ \ \ & \ \ \ \ $\Delta$ \ \ \ \ \\ \hline
  %$1$     & 0 & 0 & 0 \\ \hline
  %$2$     & 0 & 0 & 0 \\ \hline
  %$3$     & 0 & 0 & 0 \\ \hline
  %$4$     & 8 & 8 & 0 \\ \hline
  $5$     & 56 & 56 & 0 \\ \hline
  $6$     & 352 & 352 & 0 \\ \hline
  $7$     & 1760 & 1744 & 16 \\ \hline
  $8$     & 8192 & 8128 & 64 \\ \hline
  $9$     & 35712 & 35424 & 288 \\ \hline
  $10$    & 151040 & 149888 & 1152 \\ \hline
  $11$    & 624128 & 619456 & 4672 \\ \hline
  $12$    & 2547712 & 2529024 & 18688 \\ \hline
  $13$    & 10311680 & 10238848 & 72832 \\ \hline
  \end{tabular}}
\caption{\small{Comparison between the number of crossings in our
drawing and the conjectured values}}
\end{table}

\bigskip

\noindent {\bf Acknowledgements}

\noindent
The authors are deeply grateful to the anonymous reviewers for their many valuable and helpful suggestions, which have led to substantial improvements both
in the presentation of the paper and in the accuracy of the arguments. This work is supported by NSFC (grant no. 11971347, 60973014, 61303023).

\begin{figure}
\centering
\includegraphics[scale=1.0]{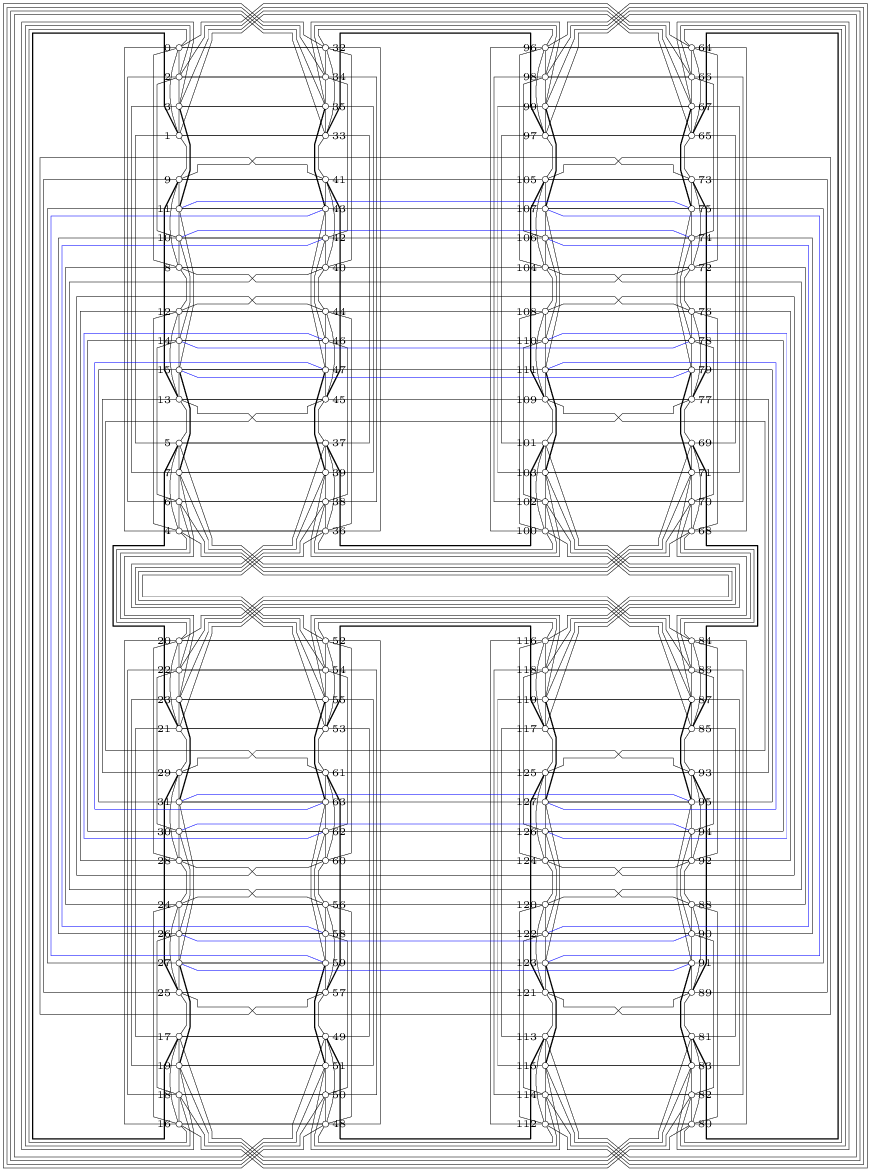}\caption{\small{The
drawing $\Gamma_7$}} \label{fig:wholedrawing7}
\end{figure}

\begin{figure}
\centering
\includegraphics[scale=1.0]{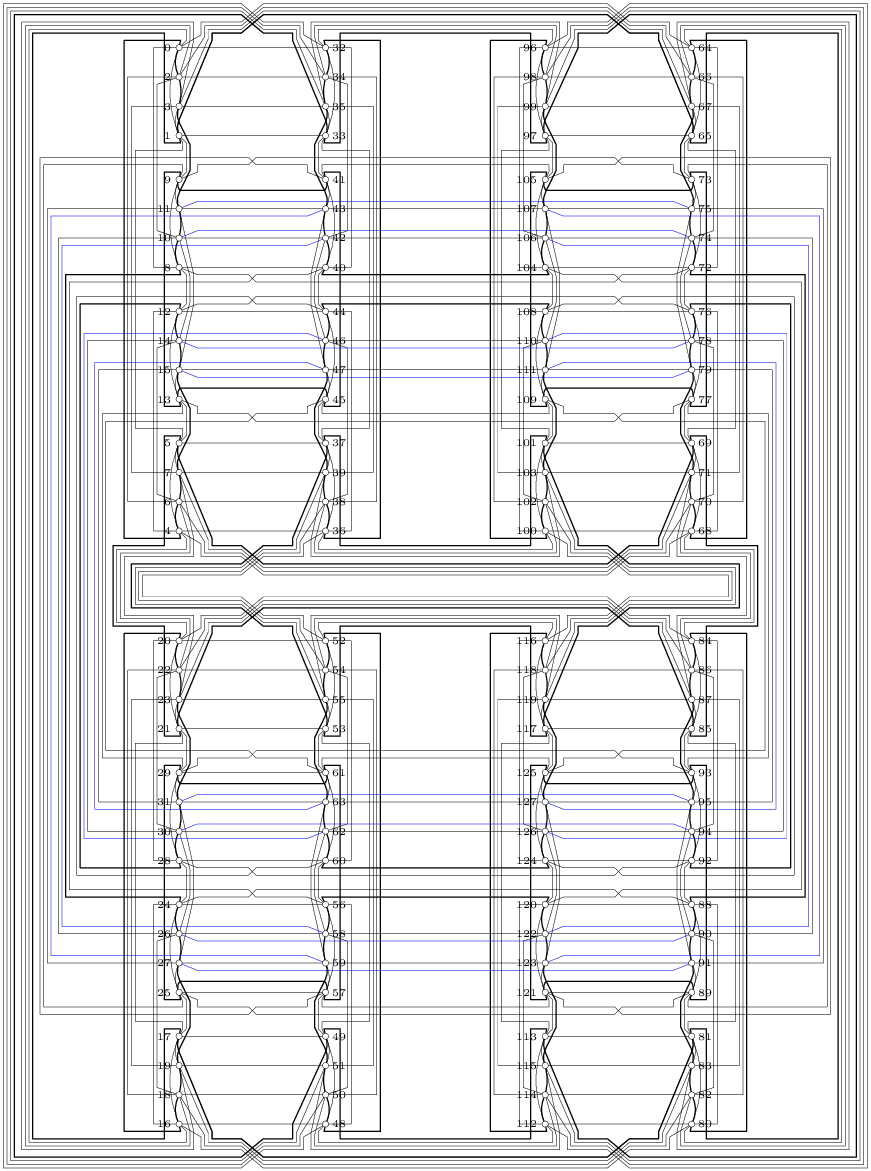}\caption{\small{The
drawing $\Gamma^{*}_7$}} \label{fig:wholedrawing77}
\end{figure}

\begin{figure}
\centering
\includegraphics[scale=1.0]{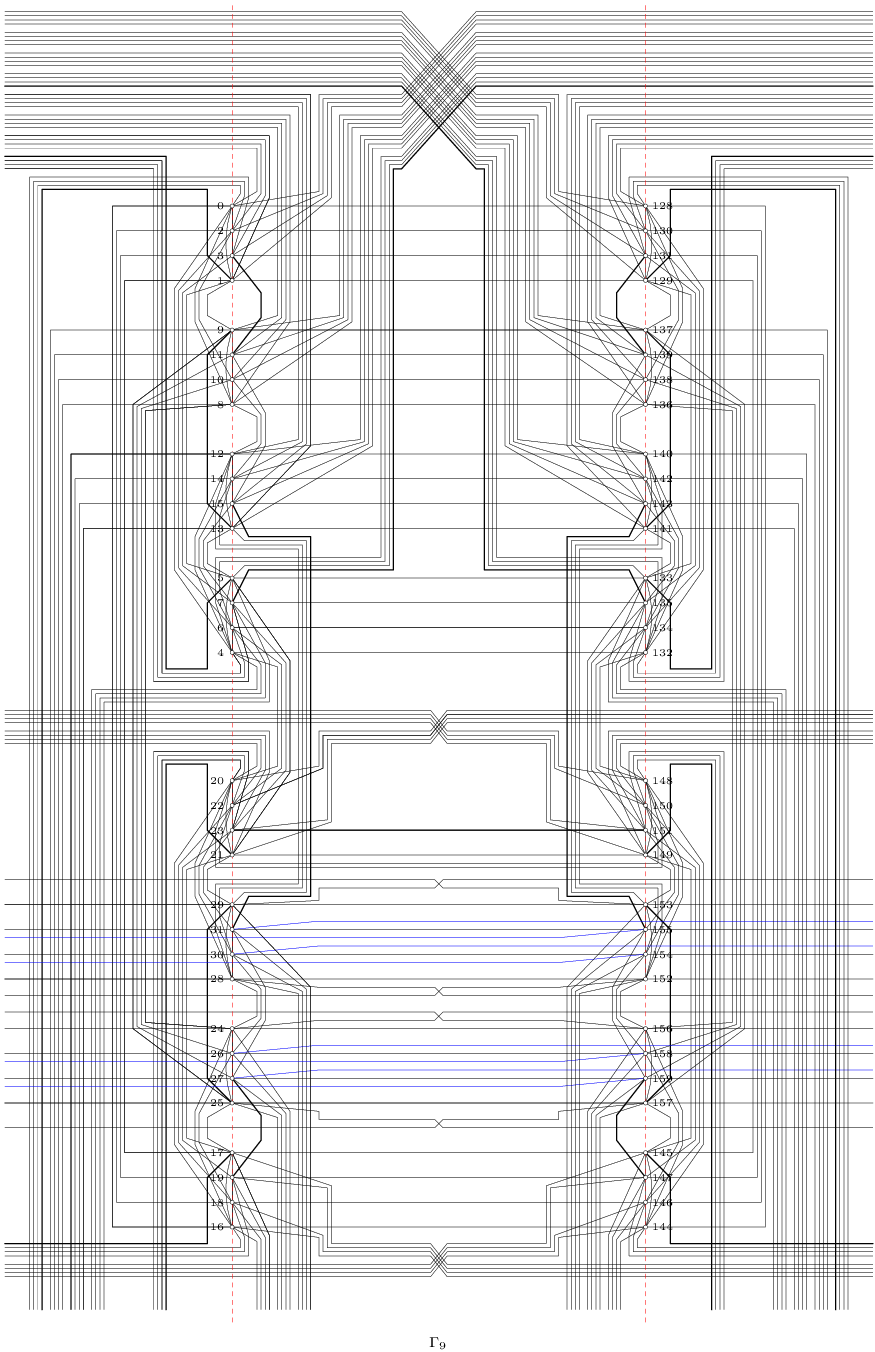}
\caption{\small{The partial drawing of $\Gamma_9$}}\label{fig:gamma
9}
\end{figure}

\begin{figure}
\centering
\includegraphics[scale=1.0]{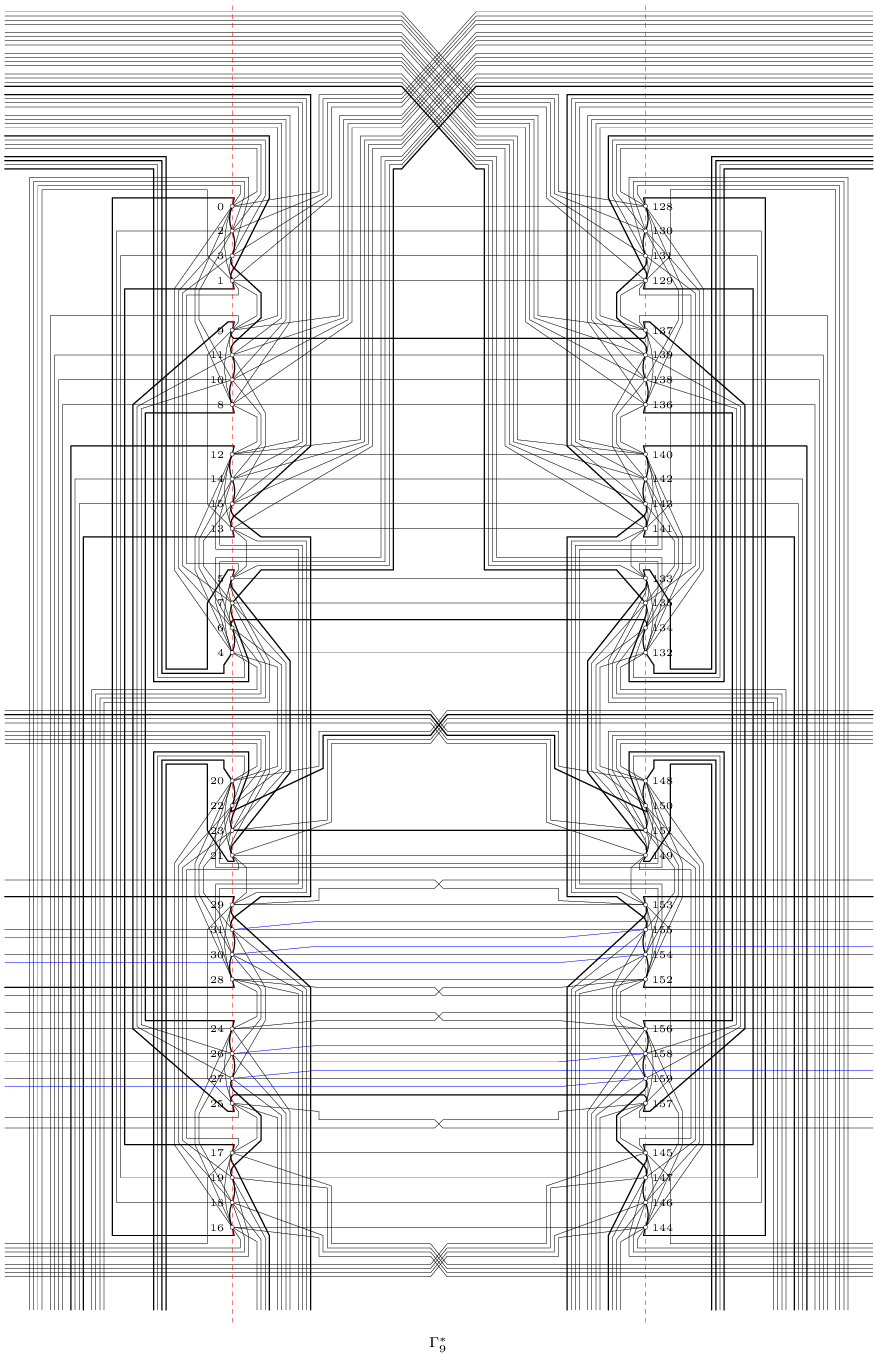}
\caption{\small{The partial drawing of
$\Gamma_9^{*}$}}\label{fig:Gamma9star}
\end{figure}

\begin{figure}
\centering
\includegraphics[scale=2.5]{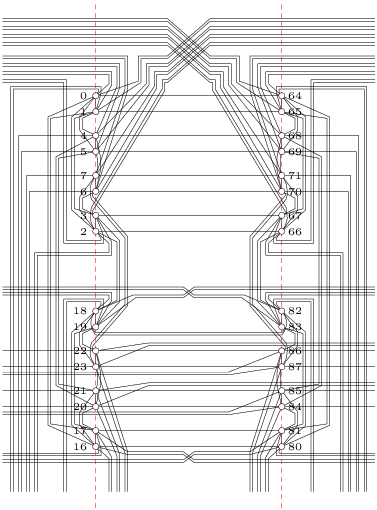}
\caption{\small{The partial drawing of
$\Gamma_8$}}\label{fig:Gamma8}
\end{figure}

\end{document}